\documentclass[11pt]{amsart}
\usepackage{enumerate,color}
\usepackage{showidx}
\usepackage{latexsym,enumerate}
\usepackage{amsmath,amsthm,amsopn,amstext,amscd,amsfonts,amssymb,mathrsfs}
\usepackage{float}
\usepackage{xcolor}
\usepackage{graphicx}
\usepackage{tikz}
\usepackage[english]{babel}
\usepackage[latin1]{inputenc}
\usepackage{times}
\usepackage[T1]{fontenc}
\usepackage{verbatim}
\usetikzlibrary{arrows,shapes}
\numberwithin{equation}{section}
\usepackage[labelformat=empty]{caption}
\def\Box{\vrule height 6pt depth 1pt width 4pt}
\date{}
\begin{document}
\begin{center}
\Large{\bf  A unified approach to symmetry for semilinear equations }
\\
\Large{\bf associated to the Laplacian 
in $\mathbb{R}^N$}
\\[0.5cm]
\end{center}
\begin{center}
{\sc A. Avila$^{1,2} $ - F. Brock$^3$}
\end{center}
\vspace*{0.5cm}
\thanks{}
\date{\today}
\small
\setcounter{footnote}{1}\footnotetext{Universidad de La Frontera, Departamento de Ingenier\'{\i}a Matem\'atica, Temuco, Chile, email: andres.avila@ufrontera.cl}
\setcounter{footnote}{2}\footnotetext{Institut f\"ur Mathematik\\Universit\"at Kassel, Germany,email: avila@mathematik.uni-kassel.de}
\setcounter{footnote}{3}
\footnotetext{University of Rostock, Institute of Mathematics, Ulmenstr.69, Haus 3, 18057 Rostock, Germany,
email: friedemann.brock@uni-rostock.de
}
\small
{\sl Abstract: }
We show radial symmetry of positive solutions to the H\'{e}non equation $-\Delta u = |x|^{-\ell} u^q $ in $\mathbb{R}^N \setminus \{ 0\} $, where $\ell \geq 0$, $q>0$ and satisfy further technical conditions.    
A new ingredient is a maximum principle for open subsets of a half space. It allows to apply the Moving Plane Method once a slow decay of  the solution at infinity has been established, that is $\lim _{|x|\to \infty } |x|^{\gamma } u(x) =L $, for some numbers $\gamma \in (0, N-2)$ and $L >0$. Moreover, some examples of non-radial solutions are given for $q> \frac{N+1}{N-3}$ and $N\geq 4$. We also establish radial symmetry for related and more general problems in $\mathbb{R}^N $ and $\mathbb{R}^N \setminus \{ 0\} $.      
\\[0.3cm]
\normalsize
\small
\renewcommand{\baselinestretch}{1.2}
\normalsize
\\[0.3cm]
\def\Box{\vrule height 6pt depth 1pt width 4pt}
\section*{1. Introduction} 
\setcounter{section}{1}
\setcounter{equation}{0}
\hspace*{0.3cm}
Due to the celebrated articles \cite{Se, GNN,GNNrn}, the Alexandroff-Serrin Moving Plane Method (MPM) has been established as a powerful tool to obtain symmetry properties of solutions to elliptic equations in bounded or unbounded  domains. See for instance the papers \cite{CafGS,BeNi,CongmingLi,LiNi1,LiNi2,terracini, DamPaRa,SeZ,Naito,DamRa,DamSci,cafflinirenberg,sciunzi}. The survey article \cite{weimingni} and the monographs \cite{Fra,ChenLi} provide further material on this subject.  
\\ 
In the case of unbounded domains in $\mathbb{R}^N $, asymptotic properties of the solution at infinity are essential to provide a starting point for the MPM.
In some papers, this problem was settled by obtaining precise estimates of the solutions involving terms of higher order, which required a lot of effort, see \cite{LiNi0,YiLi,LiNi1,Z1,Z2,Z3}. Our aim is to bypass these technicalities 
and to develop a simplified approach. At the same time, we will  obtain new symmetry results for solutions of semilinear problems associated to the Laplacian on $\mathbb{R}^N$, and in particular for homogeneous non-linearities as they appear in the 
H\'{e}non equation
\begin{equation}
\label{henon1}
-\Delta u= |x|^{-\ell } u^q,
\end{equation}
where $\ell \geq 0$ and $q>0$. 
\\
Without further ado, let us state our main results, leaving all further background information to the next section.  Throughout this paper we will use the following numbers: 
\begin{eqnarray*}
 & & q_1 (\ell)  :=  \frac{N-\ell}{N-2} ,
\qquad 
q_2 ( \ell )   :=  \frac{N+2 -2\ell}{N-2} ,
\\
& & q_{S}  :=  \left\{
\begin{array}{ll} \frac{N+1 }{N-3} & \ \mbox{ if } \ N>3
\\
 +\infty & \ \mbox{ if } \ N=3
\end{array}
\right.  ,
\\
 & & \gamma  :=  \frac{ 2-\ell }{q-1} \ \mbox{ if } \ q\not= 1, \quad \mbox{and}
\\
 & & L := \left[ \gamma (N-2-\gamma) \right] ^{1/(q-1)} .
\end{eqnarray*}
We consider solutions to the following two problems:
$$
\mbox{{\bf (P)}}  
\left\{
\begin{array}{l} 
u\in C^2 (\mathbb{R}^N \setminus \{ 0\} )\cap C(\mathbb{R}^N) , 
\\
 -\Delta u = f(|x| ,u ), \ u>0 \quad \mbox{on }\ \mathbb{R}^N \setminus \{ 0 \},
\\
 \lim _{|x|\to \infty } u(x)=0, 
\end{array}
\right.
$$
and 
$$
\mbox{{\bf (P)}$_{\mathbf{0}}$} 
\left\{ 
\begin{array}{l}
u\in C^2 (\mathbb{R}^N \setminus \{ 0\} ) , 
\\
  -\Delta u = f(|x| ,u ), \ u>0 \quad \mbox{on }\ \mathbb{R}^N
\setminus \{ 0\}  ,
\\
  \lim _{|x|\to \infty } u(x)=0, \qquad
  \lim_{x\to 0} u(x) = + \infty ,
\end{array}
\right.
$$
where  the nonlinearity $f$ satisfies 
\begin{eqnarray}
 & &  f\in C^1 \left( (0, + \infty ) \times (0, +\infty )  \right) \, ,
\label{fsmooth}
\\
 & & \mbox{the mapping $r\longmapsto f(r,u)$ is non-increasing}.
\label{fradial}
\end{eqnarray}
Note that assumption (\ref{fradial}) is indispensible when proving symmetry using the MPM or rearrangement tools. 
\\[0.1cm]
\hspace*{1cm}We will say that $u$ is {\sl radially symmetric and radially decreasing with respect to a point $x^0 \in \mathbb{R}^N $ }
if there is a function
$U\in C^1 (0, +\infty  ) $ such that $u(x ) = U(|x -x^0 |) $ for all $x\in \mathbb{R}^N \setminus \{ x^0 \} $, and $U'(r) <0 $ for $r>0 $. 
\\[0.1cm]
The main result of our paper is 
\\[0.1cm]
{\bf Theorem 1.1. } {\sl Assume that $f$ satisfies the conditions 
(\ref{fsmooth}), (\ref{fradial}) and 
\begin{eqnarray}
\label{fxluq} & & 
\left\{ 
\begin{array}{l}
f(r,u ) = r^{-\ell } u^q \left[1 + O\left( ( r^{-2} + u ^2 ) ^{\varepsilon /2 }\right) \right] ,
\\
 f_r (r,u) = - \ell r^{-\ell -1} u^q  \left[ 1 + O\left( ( r^{-2} + u ^2 ) ^{\varepsilon /2 }\right) \right] ,
\\
 f_u (r,u ) =  q r ^{-\ell} u^{q-1} \left[ 1 + O\left( ( r^{-2} + u ^2 ) ^{\varepsilon /2 }\right) \right] ,\quad \mbox{as $ r^{-2} + u ^2   \to 0 $},
\end{array}
\right.
\end{eqnarray}
where $\ell \geq 0 $, $q> 0 $, and  $\varepsilon \in (0, 1] $, and let $u$ be a solution of {\bf (P)} or 
{\bf (P)}$_{\mathbf{0}}$. Furthermore, 
assume that one of the following conditions {\bf (i)} or {\bf (ii)}
is satisfied:
\\
{\bf (i) } 
$\ \ell \in [0,2)$, $q>q_1 (\ell) $, $q\not= q_2 (\ell)$; 
\\
additionally, if $q\geq  \frac{N+2}{N-2}$, then also
\begin{equation}
\label{decay2}
 u(x) \leq c |x| ^{-\gamma } \quad \mbox{ for some $c>0$,}
\end{equation} 
and  
if $N\geq 4$, then either $\ell >0 $ and   
$q\leq q_S $, or $\ell =0 $ and $q<q_S $.
\\
{\bf (ii) } $\ \ell \in (2,N) $ and $q\in (0,q_1 (\ell ))$.
\\ 
Then 
$u$  is radially symmetric and radially decreasing w.r.t. some point $x^0 \in \mathbb{R}^N$. Moreover, if $\ell>0$, or if $u$ is a solution of {\bf (P)}$_{\mathbf{0}}$, then $x^0 =0$. }
\\[0.1cm]
The proof of Theorem 1.1 is based on the MPM. A new ingredient is a {\sl  maximum principle for open subsets of a halfspace } (see Theorem 4.2 of Section 4). It  allows to reduce the effort for asymptotic estimates at infinity. To make the MPM work, we will only need the following limit property for some number $L> 0$,
\begin{equation}
\label{simplim}
\lim_{|x|\to \infty } |x|^{\gamma } u(x) =L.
\end{equation}
Theorem 1.1 yields in particular the following result for the H\'{e}non equation.
\\[0.1cm]
{\bf Corollary 1.2. } {\sl Let $u$ be a solution of {\bf (P)} or  {\bf (P)}$_{\mathbf{0}}$, where  $f(r,u) = r^{-\ell} u^q $ and $\ell $, $q$ satisfy either {\bf (i)} or {\bf (ii)}. Then, the assertions of Theorem 1.1. hold true.} 
\\[0.1cm]
There are some situations when the conditions (\ref{fxluq}) can be considerably relaxed. See for instance \cite{GNNrn,LiNi0,LiNi1,LiNi2,YiLi,Naito}. Our next two Theorems partly overlap or extend these results. The proofs are much simpler than the proof of Theorem 1.1, and instead of the limit property (\ref{simplim}) we  merely use an {\sl upper} estimate for the solution at infinity.  First we consider the case that the solution $u$ decays faster than in (\ref{simplim}).  
\\[0.1cm]
{\bf Theorem 1.3. } {\sl  Assume that $f$ satisfies the conditions 
(\ref{fsmooth}), (\ref{fradial}) and that there are positive numbers $u_0 , \, r_0, \, d_1, \, d_2 $ such that  
\begin{eqnarray}
\label{fxluqsimple} & & 
\left\{ 
\begin{array}{l}
0< f(r,u ) \leq d_1  r^{-\ell } u^q 
\\
 f_u (r,u ) \leq   d_2  r ^{-\ell} u^{q-1} 
\end{array}
\right. 
, \quad \mbox{for $\ 0<u<u_0$ and $r>r_0 $,}
\end{eqnarray}
where
\\
{\bf (i') } $\ \ell \in [0,2)$ and $q >q_1 (\ell)$. 
\\
Furthermore, let $u$ be a solution of   {\bf (P)} or  {\bf (P)}$_{\mathbf{0}}$  with
\begin{equation}
\label{faster}
\lim_{|x|\to \infty } |x|^{\gamma } u(x) = 0.
\end{equation}
Then the assertions of Theorem 1.1 hold true.} 
\\[0.1cm]
Next we consider cases where the parameters $\ell $ and $q$ fall out of the ranges {\bf (i')} or {\bf (ii)}. As we shall see in  Section 3,  the solutions show fast decay in the cases {\bf (iv)}--{\bf (vi)} below. Note however the curious fact that we do not need any asymptotic estimate for the solution in case {\bf (iii)} (see the proof of Theorem 1.4 in Section 5).
\\[0.1cm]   
{\bf Theorem 1.4. } {\sl Assume that $f$ satisfies the conditions 
(\ref{fsmooth}), (\ref{fradial}) and (\ref{fxluqsimple}), and 
that one of the following conditions is satisfied:
\\
{\bf (iii) } $\ \ell =2 $ and $q>1 $; 
\\
{\bf (iv)} $\ \ell \in (2,N) $ and $q= q_1 (\ell )$; 
\\
{\bf (v)} $\ \ell \in (2,N) $ and $q> q_1 (\ell )$; or
\\
{\bf (vi)} $\ \ell \geq N$ and $q>0$.
\\ 
Then, if $u$ is a solution of problem {\bf (P)} or of {\bf (P)}$_{\mathbf{0}}$, it is radially symmetric and radially decreasing w.r.t.  $0$.
}
\\[0.1cm]
{\bf Remark 1.5. (a) } 
 Our results cover a large range of values $(\ell,q)$. Indeed, assume that there are positive numbers $u_0 $, $r_0 $ and $d_3$ such that  the nonlinearity $f$ satisfies
\begin{equation}
\label{fsub}
f(r,u)\geq d_3 r^{-\ell} u^q \quad \mbox{for $0<u<u_0 $ and $r>r_0$.}
\end{equation} 
Then, problems {\bf (P)} and {\bf (P)}$_{\mathbf 0}$ do not have positive solutions when 
\\
{\bf (vii) } $\ \ell\in [0,2)$ and $q\in (0, q_1 (\ell)] $, see \cite{BidautPoho}, Theorem 3.3 (ii) and Theorem 3.4 (ii);
\\
{\bf (viii) } $\ \ell=2$ and $q\in (0,1)$, see \cite{BidautPoho}, Theorem 3.4 (ii).   
\\ 
{\bf (b)} We will give examples of {\sl non-radial} solutions to problem {\bf (P)$_{\mathbf 0}$} when $\ell \in [0,2)$, $N\geq 4$ and $q>q_S $ in Section 6. 
Therefore, it is not clear  under which conditions the solutions are radial in these cases. 
It would be also very interesting to find non-symmetric solutions for problem {\bf (P)}  when $\ell\in [0,2)$, $q>q_S $ and $N\geq 4$. 
The figure below illustrates the different values  
in the $(\ell ,q)$-plane for $N\geq 4$ in our Theorems 1.1-- 1.3. The grey trapezoid is the region of nonexistence.  

\begin{figure}[ht]
\centering
\includegraphics[scale=0.45]{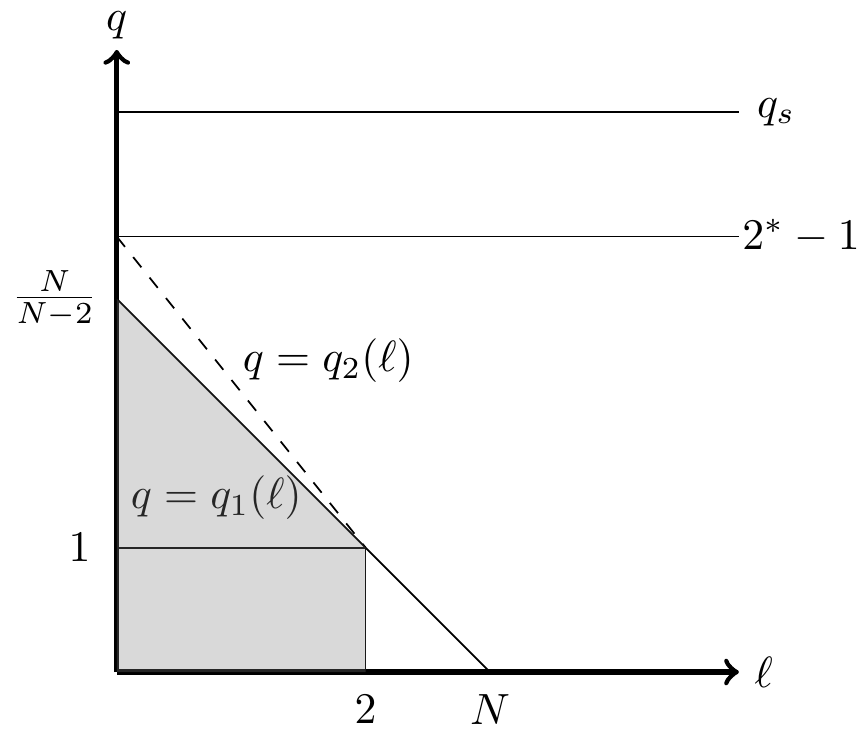}
\end{figure} 

Now we outline the content of the paper.
In section 2, we collect some preliminary material and results that are available in the literature for problems {\bf (P)} and {\bf (P)}$_{\mathbf 0} $. In section 3, we obtain lower and upper bounds at infinity for the solutions of our problems. Furthermore, we show that the solutions satisfy the limit property (\ref{simplim}) under assumptions of Theorem 1.1.  In section 4, we prove Theorem 4.2, which is a maximum principle for open subsets of a half space. Section 5 deals with the proofs of Theorems 1.1, 1.3 and 1.4, which are based on the MPM and the results of sections 3 and 4. Finally, in section 6 we give some examples of non-radial solutions to problem {\bf (P)}$_{\mathbf 0}$ when $q>q_S $.  
\section*{2. Preliminaries}
\setcounter{section}{2}
\setcounter{equation}{0}
Our work was inspired by articles of    
H. Zou related to the Lane-Emden equation, see \cite{Z1,Z2,Z3}, and some progress concerning 
 asymptotic estimates of positive solutions to elliptic equations in exterior domains. See \cite{BVV,BidautPoho,MitPoh,PoQuSo,AvilaBrock,PhanSouplet} and the references cited therein. 

H. Zou  considered the following problem:   
\begin{equation}
\label{Zouproblem}
\left\{
\begin{array}{l} 
 u\in C^2 (\mathbb{R}^N ),
\\
 -\Delta u= f(u), \quad u>0 \quad \mbox{on }\ \mathbb{R}^N,
\\
 \lim_{|x|\to \infty} u(x)=0,
\end{array}
\right.
\end{equation}
where $f$ is smooth with $f(u)\sim u^q $ near $u=0$ for some $q>1$. A key result is
\\[0.1cm]
{\bf Theorem A. } {\sl (see \cite{Z3}, Theorem 1.2 and Lemma 2.1)
 Let $u$ be a solution of  problem (\ref{Zouproblem}), 
where $N\geq 3 $, 
$f\in C^1([0,+\infty))$, $f(u)>0$ for $u>0$, and    
\begin{equation}
\label{fasymp}
 f(u) = u^q + O(u^{q+\varepsilon} ) , \quad f^{\prime} (u) = qu^{q-1} + O( u^{q-1+\varepsilon } ),  \quad \mbox{for  $ u\in (0, u_0 )$,}
\end{equation}
for some positive numbers $\varepsilon $ and $u_0 $ and
$$
 1<q< \frac{N+2}{N-2}.
$$
Then $u$ is radially symmetric about some point.}
\\[0.1cm] 
It follows from the proof in \cite{Z3} that Theorem A also holds  without the positivity assumption for $f$.   Furthermore, solutions are radially symmetric as well for the range 
$$
\frac{N+2}{N-2} <q< \left\{ 
\begin{array}{ll} +\infty & \mbox{ if $N=3$}
\\
\frac{N+1}{N-3} & \mbox{ if $N\geq 4$} 
\end{array}
\right.
,
$$ 
provided that they satisfy the following estimate from above at infinity,
\begin{equation}
\label{asympuq}
u(x)\leq C|x|^{-2/(q-1)} \, , \quad |x|\geq 1 ,
\end{equation}
for some $C>0$. In the special case $f(u)=u^q $, this estimate was stated in  Theorem 1.1 in \cite{Z1}, but the arguments used in \cite{Z1,Z2,Z3} also carry over to the general case. H. Zou also showed (\ref{asympuq}) under additional conditions on the solution, see Theorem 1.2 in \cite{Z2}. Furthermore, we mention that Z. Guo \cite{Guo} extended the results of \cite{Z1,Z2,Z3} to $q\geq \frac{N+1}{N-3}$, $N\geq 4$. More precisely, he showed that all positive $C^2 $--solutions of the Lane--Emden equation $-\Delta u =u^q $ in $\mathbb{R}^N $ are radially symmetric in the following two cases: 
\\
{\bf (a) } $\ N\geq 5$, $q\geq \frac{N}{N-4}$ and  $u$ satisfies \begin{equation}
\label{precasymp}
\lim_{|x|\to \infty } |x|^{2/(q-1) } u(x) = \lambda ,
\end{equation}
where 
\begin{equation}
\label{lamb}
\lambda = \left[ \frac{2}{q-1} \left( N-2 -\frac{2}{q-1} \right) \right] ^{1/(q-1)}.
\end{equation}
{\bf (b) } $\ N\geq 4$, $\frac{N+1}{N-3}\leq q<\frac{N}{N-4}$ and $u$ satisfies  (\ref{precasymp}) and
\begin{equation}
\label{morepreciseu}
\lim_{|x|\to \infty } |x|^{1-(\mu +N)/2 } \left( |x|^{2/(q-1)}u(x)-\lambda \right) =0,
\end{equation}
where $\lambda $ is given by (\ref{lamb}) and  
$$
\mu = \frac{4}{q-1}+4 -2N.
$$ 
\hspace*{0.3cm}It is natural to ask about qualitative properties of solutions when the right-hand side $u^q $ is replaced by a more general term which is homogeneous in $u$ and $|x|$.  
The resulting PDE is  the so-called {\sl H\'{e}non equation},
\begin{equation}
\label{henon}
-\Delta u = |x|^{-\ell} u^q , 
\end{equation} 
where 
$\ell \in \mathbb{R}$ and $q >0$, and 
it  appears in Geometry and Physics. It also serves as a model for many other semilinear problems, and it has been extensively studied, both in bounded and unbounded domains,   see e.g. \cite{GidSpruck,BVV,BidautPoho,DaDuGuo,Dengetal, PhanSouplet}.
\\ 
\hspace*{0.3cm}Let $\Omega $ be a domain in $\mathbb{R}^N$ with $0\not\in \Omega $, and consider the problem
\begin{equation}
\label{henon1}
\left\{ 
\begin{array}{l}
u\in C^2 (\Omega ),
\\
-\Delta u = |x|^{-\ell} u^q ,
\ \ u>0 \quad \mbox{in $ \Omega $.}
\end{array}
\right.
\end{equation} 
For the asymptotic properties of the solutions  near $0$ and infinity, 
 numbers $q_1 (\ell )$, $q_2 (\ell )$, $q_S $, $\gamma $ and $L$ defined in Section 1 play an important role. The following result refers to the subcritical case.
\\[0.1cm]
{\bf Theorem B. }{\sl (see \cite{GidSpruck}, Theorem 3.4) Let $u$ be a solution of (\ref{henon1}), $\ell <2$, $q_1 (\ell ) <q <\frac{N+2}{N-2}$ and $q\not= q_2 (\ell )$.
\\
If $\Omega = B_1 \setminus \{ 0 \} $, then either $x=0$ is a removable singularity of $u$, or $x=0$ is a non-removable singularity and 
\begin{equation}
\label{asymp1}
\lim_{x\to 0} |x| ^{\gamma } u(x)= L . 
\end{equation}
On the other hand, if $\Omega = \mathbb{R}^N \setminus \overline{B_1 }$, then either
\begin{equation}
\label{asymp2}
\lim_{|x|\to \infty } |x| ^{N-2 } u(x)= \lambda , \quad \mbox{(fast decay),}
\end{equation}
for some constant $\lambda >0$, or
\begin{equation} 
\label{asymp3}
\lim_{|x|\to \infty } |x| ^{\gamma } u(x)= L , \quad \mbox{ (slow decay).}
\end{equation}
}
\hspace*{0.3cm}There is also a result in the  supercritical case $q>\frac{N+2}{N-2}$. 
\\[0.1cm]
{\bf Theorem C. } {\sl (see \cite{BVV}, Theorems 3.2, 3.3 and Remark 3.2) Let $u$ be a solution of (\ref{henon1}). 
If $\Omega = B_1 \setminus \{0\} $, $\ell <2$ and $q> 
\max \{ q_1 (\ell),\,\frac{N+2}{N-2} \} $, and if 
\begin{equation}
\label{uless}
u(x)\leq c|x|^{-\gamma }, \quad \mbox{for some $c>0$,
}
\end{equation} 
then either $0$ is a removable singularity of $u$ or 
\begin{equation}
\label{asymp4}
\lim_{x\to 0} |x|^{\gamma } u(x) = V(\theta ), \quad 
\mbox{uniformly in $\theta = \frac{x}{|x|}$, }
\end{equation}
where $V\in C^2 (\mathbb{S}^{N-1})$ is a positive solution of
\begin{equation}
\label{spherical}
-\Delta _{\theta } V + \gamma (N-2-\gamma ) V= V^q \quad \mbox{on $\mathbb{S}^{N-1} $,}
\end{equation}
($\Delta _{\theta } $= Laplace-Beltrami operator on the sphere 
$\mathbb{S}^{N-1} $). 
\\
Finally, if $\Omega = \mathbb{R}^N \setminus \overline{B_1 }$, and if $\ell$ and $q$ are as above and $u$ satisfies (\ref{uless}), then either 
\begin{equation}
\label{asymp5}
\lim_{|x|\to \infty } |x|^{N-2 } u(x) =\lambda \, , \qquad \qquad \qquad \qquad \qquad \mbox{(fast decay)},
\end{equation}
 or
\begin{equation}
\label{asymp6}
\lim_{|x|\to \infty } |x|^{\gamma } u(x) = V(\theta )\, , \ \ \ \mbox{uniformly in $\theta = \frac{x}{|x|} \, ,  \quad $ (slow decay).}
\end{equation} 
}
Note that an easy application of \cite{BVV}, Theorem 6.1, shows that, if $q<q_S $, then (\ref{spherical}) has only the constant solution $V(\theta ) \equiv L $, (compare also Lemma 3.10 below). However, it is unclear under which conditions the estimate  
(\ref{uless}) holds in the case $q>\frac{N+2}{N-2}$. 
\\
The significance of the number $q_2 (\ell)$ can be best understood by  looking at  {\sl radial } solutions for the problem (\ref{henon1})  when 
$$
\Omega = \mathbb{R}^N \setminus \{ 0 \} .
$$
Assume again $\ell<2 $ and $q>q_1 (\ell )$. The radial solutions have been classified in \cite{GidSpruck}, Appendix A:
\\ 
{\bf 1.} The basic solutions  
\begin{equation}
\label{basic1}
u(x)= L |x|^{-\gamma };
\end{equation}
{\bf 2.} A one-parameter family of solutions in the case 
$q\in (q_1 (\ell ), q_2 (\ell))$ which satisfy (\ref{asymp1})  and (\ref{asymp2});
\\ 
{\bf 3.} Another one-parameter family of solutions in the case $q\in (q_2 (\ell ), +\infty )$  which are $C^2 $ and satisfy
(\ref{asymp3}); 
\\
{\bf 4.} Two further types of radial solutions in the critical case
$
q= q_2 (\ell)  $: 
\\
the fast-decay solution 
\begin{equation}
\label{critsol1}
u(x) = \left( \frac{\mu \sqrt{(N-\ell)(N-2)}}{ \mu ^2 + |x|^{2-\ell}}
\right) ^{(N-2)/(2-\ell)} \quad (\mu >0),
\end{equation} 
and the slow-decay solution 
\begin{equation}
\label{critsol2}
u(x) = |x|^{(2-N)/ 2}  \Psi (\ln |x|),
\end{equation}
where $\Psi $ is a strictly positive and periodic  solution of the ODE
$$
\Psi ^{\prime \prime } (t) -\left( \frac{N-2}{2} \right)^2 \Psi (t) + \left( \Psi (t)\right)  ^{(N+2-2\ell)/(N-2)} =0 \quad (t\in \mathbb{R}),
$$ 
which oscillates about the value
$
\Psi _0 = \left( \frac{N-2}{2} \right) ^{(N-2)/(2-\ell)} .
$
\\
Note that if $\ell <2$ and $q\in [0,q_1 (\ell)  ]$, then  problem (\ref{henon1}) with $\Omega = \mathbb{R} ^N \setminus \{ 0\} $ has no solution.
\\
\hspace*{0.3cm} We can see that  the asymptotic behavior at $0$ and infinity of arbitrary solutions of (\ref {henon1}) can be read off from the asymptotics  of the {\sl radial} ones. Moreover, if 
$\ell \geq 0$, then the right-hand side of (\ref{henon}) has the right monotonicity behavior for a successful application of the MPM. 
\\[0.1cm]
{\bf Remark 2.1. 
 (a)} Apart from the H\'{e}non equation, there are further  well-known PDE  covered by our results. Two examples are:
\\
{\bf 1.} The generalized Matukuma equation
\begin{equation}
\label{matuk}
-\Delta u = \frac{|x|^{\lambda -2}}{(1+|x|^2)^{\lambda /2} }u^q , \quad (\lambda >0),
\end{equation}
see \cite{YiLi,LiNi0} and Theorem 1.3, {\bf (iii)};
\\
{\bf 2.} the scalar curvature equation
\begin{equation}
\label{scalcurv}
-\Delta u= K(|x|) u^{\frac{N+2}{N-2}}, 
\end{equation}
when $
K(r) \sim r^{-\ell} $ for large $r$, $(\ell >0$),
see \cite{LiNi1} and Theorem 1.1, {\bf (i)} and 1.3, {\bf (v)} and {\bf (vi)}.
\\
{\bf (b)} We have excluded the borderline case
$$
l\in [0,2) \ \mbox{ and }\ q= q_2 (l) ,
$$
in Theorem 1.1, 
because the proof of radial symmetry would require further tools (compare with Remark 3.8). Note that the special case of the equation
\begin{equation}
\label{CKN}
-\Delta u= |x|^{-\ell} u^{q_2 (\ell)} 
\end{equation}
has received a lot of attention during the last years, since it  is related to some Caffarelli-Kohn-Nirenberg inequalities, see \cite{Hsia} and the references cited therein. Combining the MPM and appropriate Kelvin transformations, one can prove that all solutions of Problem  {\bf (P)}$_{\mathbf 0} $ for equation (\ref{CKN}) are radially symmetric and radially decreasing, see Theorem 1.1 in \cite{Hsia}. However, this method seems not  to be applicable under the general assumptions (\ref{fxluq}). 
\\[0.1cm]
{\bf Remark 2.2.} {\bf (a)} Y. Naito \cite{Naito} studied problem {\bf (P)} when $f(r,u)= \varphi (r) g(u) $ with continuous, non-increasing and non-negative $\varphi $ satisfying $\varphi (r)= O(r^{-\ell})$ as $r\to \infty $, for some $\ell \in [0,2] $, and $g\in C^1 [0,+\infty )$ satisfying $g'(u) = O (u^{q-1} ) $ as $u\to 0$, for some $q>(N-\ell )/(N-2)$. Using the MPM and a maximum principle in unbounded domains (compare Lemma 4.1 of Section 4 below) he proved radial symmetry of $u$ provided that it  decays fast at infinity, that is,
$$
u(x) = \left\{ 
\begin{array}{ll}
o(|x| ^{-\gamma } ) & \mbox{ as $|x|\to \infty ,\ $ if $\  \ell \in [0, 2)$} 
\\
 o ((\log |x|) ^{-1/(q-1)} ) & \mbox{ as $|x|\to \infty , \ $ if $\ \ell =2 $}
\end{array}
\right.
,
$$
see Theorem 1 and Corollary 2 of \cite{Naito}. However, he did not consider the more difficult case  that $u$ has the slow decay (\ref{simplim}). This fact was an important motivation for our study. 
\\
{\bf (b)} An alternative approach to symmetry is based on the MPM together with integral estimates, see  \cite{terracini,cafflinirenberg,sciunzi}. But this method has been applied only to situations when the solution belongs to the space $L^{2N/(N-2)} (\mathbb{R}^N \setminus \overline{B_1 })$.  
\\[0.1cm]
\hspace*{0.3cm}Throughout this paper, we will use the following notation:
\begin{eqnarray*}
 & & \mathbb{R}^N \ni x= (x_1 , x_2 , \ldots , x_N ) \equiv (x_1 , x'),
\\
 & &  \mathbb{R}^N _+ := 
\{ x=(x_1 , x^{\prime} )\in \mathbb{R}^N , \, x_1 >0 \} ,
\\
 & & B_R := \{ x\in \mathbb{R}^N :\, |x|<R\}, \quad
B_R ^C := \mathbb{R}^N \setminus \overline{B_R }, \quad R>0,
\\ 
 & & u_+ := \max \{ u, 0\} , \ \mbox{ for any function or number u}.
\end{eqnarray*}
\section*{3. Asymptotic estimates at infinity }
\setcounter{section}{3}
\setcounter{equation}{0}
In this section, we obtain some estimates for the solutions of our problems at infinity. 
First we study a related problem in the exterior of a ball:
\begin{equation}
{\bf (Q)}
\left\{ 
\begin{array}{l}
u\in W^{2,N} _{loc}  (B_R ^C )\cap C\left( \overline{ B_R ^C}  \right) ,
\\
-\Delta u = f(x,u), \ \ u> 0 \quad \mbox{on } B_R ^C,
\\
\lim_{|x|\to \infty } u(x)=0,
\end{array}
\right.
\end{equation}
where $f : B_R ^C \times [0,+\infty ) \to [0, \infty )$ is measurable, locally bounded, and continuous in the
second variable. 
We begin with 
some estimates from below and from above for solutions of problem {\bf (Q)}.   
\\[0.1cm]
{\bf Lemma 3.1. } {\sl  There exists a number $c_1 >0$ such that} 
\begin{equation} 
\label{frombelow}
u(x)\geq c_1 |x| ^{2-N} \quad \mbox{on }\ B_{2R} ^C .
\end{equation}
{\sl Proof: } Since  $f\geq 0$ and 
$
\Delta (|x|^{2-N}) =0$ ,
the assertion follows from a simple comparison argument.
$\hfill \Box $
\\[0.1cm] 
{\bf Lemma 3.2. } {\sl (see \cite{AvilaBrock}, Theorem 8 ) $ \, $ Let $\, \ell \in [0,2)$, $q\in \big(q_1 (\ell ), \frac{N+2}{N-2}\big) $ and assume that 
\begin{equation}
\label{thm8}
\lim_{|x| \to \infty ,\, u\to 0 } |x|^{\ell } u^{-q} f(x,u)= \mu 
\end{equation}
for some number $\mu >0$. 
Then, there is a number $c_2 >0$ such that 
\begin{equation} 
\label{decaycrucial}
u(x) \leq c_2 |x|^{-\gamma } \quad \mbox{ on $B_{2R} ^C$.}
\end{equation}
}
{\bf Lemma 3.3. } {\sl Let $\ell \in (2,N)$ and $q\in (0,q_1 (\ell))$
and assume that
\begin{equation}
\label{l2-1}
f(x,u) \leq d_1 |x|^{-\ell } u^q \quad \mbox{for $u>0$ and $|x|\geq R$.}
\end{equation}
Then (\ref{decaycrucial}) holds.  
Moreover, there exists a number $c_3 >0$ such that 
\begin{equation}
\label{lowerboundii}
u(x)\geq c_3 |x|^{-\gamma } \quad \mbox{ on $B_{2R} ^C$.}
\end{equation}
}
{\sl Proof: } (\ref{decaycrucial})  was proved 
in \cite{AvilaBrock}, Theorem 4 and (\ref{lowerboundii}) was shown in \cite{BidautPoho}, Theorem  3.1. 
$\hfill \Box$
\\[0.1cm]
{\bf Lemma 3.4.} {\sl (see \cite{AvilaBrock}, Theorem 4) $\, $ Assume that $f$ satisfies (\ref{l2-1}) and that either
{\bf (v)} or 
{\bf (vi)} holds.
Then, there is  $c_4 >0$ such that}
\begin{equation}
\label{decaystrong}
u(x)\leq c_4 |x|^{2-N}   \quad \mbox{ on $B_{2R} ^C$.}
\end{equation}  
{\bf Lemma 3.5. } {\sl Assume that 
condition 
{\bf (iv)} holds.
Then, there is  $c_5 >0$ such that 
\begin{equation}
\label{frombelow}
u(x) \geq c_5 |x|^{2-N} \log |x| \quad \forall x\in \mathbb{R}^N \setminus \overline{B_2 } .
\end{equation}
Moreover, 
for every $\varepsilon \in (0, N-2)$ there exists a number $c(\varepsilon )>0$ such that}
\begin{equation}
\label{decayalmoststrong}
u(x) \leq c(\varepsilon ) |x| ^{\varepsilon +2-N}   \quad \mbox{ on $B_{2R} ^C$.}
\end{equation}
{\sl Proof: } 
In view of Lemma 3.1 and since $q= q_1 (l) = \frac{N-\ell}{N-2}<1$, we have for all $x\in \mathbb{R}^N \setminus \overline{B_1} $,    
\begin{eqnarray*}
-\Delta u & \geq  & C_1 |x|^{-\ell } u^{q} =  C_1 |x|^{-\ell } u^{\frac{N-\ell}{N-2}}
\\
 & \geq & C_2 |x|^{-\ell} |x|^{\ell -N} =  C_2 |x|^{ -N} ,
\end{eqnarray*} 
for some positive numbers $C_1$ and $C_2 $. 
Now (\ref{frombelow}) follows from \cite{BidautPoho}, Proposition 2.7. The second assertion (\ref{decayalmoststrong}) was proved in
\cite{AvilaBrock}, Theorem 4.
$\hfill \Box $
 
\vspace{0.3cm}

Next we want to show that solutions 
 $u$ of {\bf (P)} and of {\bf (P)}$_{\mathbf{0}}$ satisfy 
\begin{equation}
\label{liminfinity}
\lim_{|x|\to +\infty } |x|^{\gamma } u(x) = L , \quad \mbox{for some $L\geq 0$,}
\end{equation}
under condition (\ref{fxluq}) and either {\bf (i)} or {\bf (ii)}.
Our proofs depend on the previous asymptotic estimates of the solutions and on ideas of the papers \cite{Z1,Z2,Z3}. 
\\[0.1cm]
{\bf Lemma 3.6. } {\sl Assume that $f$ satisfies (\ref{fxluq}) and that either {\bf (i)} or {\bf (ii}) holds. Then, there is a constant $c_6 >0$ such that
\begin{equation}
\label{ineq1}
|D^{\alpha }  u (x)|  \leq  c_6  |x| ^{-\gamma -|\alpha |}
\qquad \forall x\in \mathbb{R}^N \setminus \overline{B_1 }, 
\end{equation}
 for every multi--index $\alpha $ with $|\alpha | \leq 3$.}
\\[0.1cm]
{\sl Proof: } By Lemma 3.2
we have that  
$$
|\Delta u (x) | \leq C_1 |x| ^{-\gamma -2} \quad \mbox{for $x\in B_1 ^C $,}
$$ 
for some positive constant $C_1 $.
Then standard elliptic estimates show (\ref{ineq1}) for every multiindex $\alpha $ with $|\alpha |\leq 2 $.
Next, differentiating the PDE  gives \\
$-\Delta u_{x_i } = f_r (|x|,u ) x_i |x| ^{-1} + f_u (|x|,u ) u_{x_i } $,
and applying once more Lemma 3.2, one obtains
$$
|\Delta u_{x_i } | \leq C_2  |x| ^{-\gamma -3 } \quad \mbox{for $x\in B_1 ^C$, $i=1,\ldots ,N $,}
$$
 for some constant $C_2 >0$. This leads
to (\ref{ineq1}) for every multi--index $\alpha $ with $|\alpha | =3 $.
$\hfill \Box $

\vspace{0.1cm}

In the following, let
\begin{equation}
\label{defv}
v(x) := |x|^{\gamma } u(x)  \quad (x\in \mathbb{R}^N ).
\end{equation}
We will also write  $v(r, \theta ) =v(x)$, where $(r, \theta )$ are spherical coordinates,
$\theta \in \mathbb{S} ^{N-1} $, $r=|x| $, and $\mathbb{S} ^{N-1} $ is the $(N-1)$--unit sphere.
It is then clear that estimates (\ref{ineq1}) imply that
\begin{equation}
\label{ineq3}
|D^{k_1 } _r D^{k_2 } _{\theta } v|\leq C r ^{-k_1 }  \quad \forall r>0 , \ \theta \in \mathbb{S}^{N-1} , 
\end{equation}
 for some $C>0 $ and 
 for any two non--negative integers  $k_1 $ and $k_2 $  with $k_1 + k_2 \leq 3 $.
\\[0.1cm]
{\bf Lemma 3.7. } {\sl Assume that either {\bf (i)} or {\bf (ii)} holds.
Then, we have
\begin{eqnarray}
\label{limv1}
 & & \lim_{r\to \infty } r v^{\prime} (r, \theta ) =0,
 \\
 \label{limv2}
 & & \lim_{r\to \infty } r^2 v^{\prime \prime} (r, \theta ) =0,
\end{eqnarray}
where $v^{\prime}$ is the derivative of $v$ w.r.t. the radius $r$.
Furthermore, the convergence is uniform in 
$C^{\tau } (\mathbb{S} ^{N-1} )$ for any
$\tau \in (0,1 )$.}
\\[0.2cm]
{\sl Proof :} We proceed in two steps, similarly as in \cite{Z1} and \cite{Z3}.
\\
{\sl Step 1 :} We claim:
\begin{eqnarray}
\label{estv1}
 & &
 r \int_{\mathbb{S} ^{N-1} }  (v^{\prime} (r, \theta ) ) ^2 \, d\theta \ \in \ L^1 (0, + \infty ) ,
 \\
\label{estv2}
 & & r^3 \int_{\mathbb{S} ^{N-1} }  (v^{\prime \prime} (r, \theta ) ) ^2 \, d\theta \ \in \ L^1 (0, + \infty ).
\end{eqnarray}
We first show (\ref{estv1}).
By direct calculation, we find that
\begin{equation}
\label{sfer1}
- v^{\prime \prime} - (N-1-2\gamma ) \frac{v^{\prime}}{r} -
\frac{\Delta _{\theta } v }{ r^2 } +
\gamma (N-2-\gamma ) \frac{v}{r^2 } =
r^{\gamma } f(r, r ^{-\gamma } v) ,
\end{equation}
where $\Delta _{\theta } $ denotes the Laplace--Beltrami operator on $\mathbb{S} ^{N-1} $.
Multiplying (\ref{sfer1}) by $r ^2 v^{\prime} $ and integrating over $(0,R) \times \mathbb{S} ^{N-1}$, $R>0$, gives
\begin{eqnarray}
\nonumber
J & := & \int_0 ^R \int_{\mathbb{S} ^{N-1} } r ^{\gamma +2 } f(r, r ^{-\gamma }v ) v^{\prime}
\\
\nonumber
 & = &
- \int_0 ^R \int_{\mathbb{S} ^{N-1} } v^{\prime \prime} v^{\prime} r^2
- (N-1-2\gamma )\int_0 ^R \int_{\mathbb{S} ^{N-1} } r(v^{\prime} )^2
\\
\label{rtheta1}
 & & + \gamma (N-2-\gamma ) \int_0 ^R \int_{\mathbb{S} ^{N-1} }  v^{\prime} v - \int_0 ^R \int_{\mathbb{S} ^{N-1} }  v^{\prime} \Delta _{\theta } v.
\end{eqnarray}
Using integration by parts, we find that the right-hand side  of (\ref{rtheta1}) 
equals to
\begin{eqnarray*}
 & &  - (N-2 -2\gamma )\int_0 ^R \int_{\mathbb{S} ^{N-1} } r(v^{\prime}) ^2
- \int_{\mathbb{S} ^{N-1} } \frac{r^2 (v^{\prime}) ^2 }{2} \Big| _0 ^R
\\
 & & + \gamma (N-2-\gamma ) \int _{\mathbb{S} ^{N-1} } \frac{v^2 }{2} \Big| _0 ^R +   \frac{1}{2} \int_{\mathbb{S} ^{N-1} }
|\nabla _{\theta } v  | ^2 \Big| _0 ^R
\\
 & =: &  - (N-2 -2\gamma )\int_0 ^R \int_{\mathbb{S} ^{N-1} } r(v^{\prime}) ^2  +I_1 + I_2 + I_3 .
\end{eqnarray*}
By (\ref{ineq3}) we find that the integrals $I_i $, ($i=1,2,3$), are uniformly bounded in $R$. Since by assumption $N-2-2 \gamma \not= 0$, it
remains to show that the left-hand side $J$ in (\ref{rtheta1}) is uniformly bounded in $R$. Using our assumptions on $f$ and (\ref{ineq3}), we have
\begin{eqnarray}
\nonumber
|J| & \leq  & C_1 + \left| \int_1 ^R \int_{\mathbb{S} ^{N-1} } v^q v^{\prime} \right| + \left|
\int_1 ^R \int_{\mathbb{S} ^{N-1} } v^q O( (r^{-2} + r ^{-2\gamma } v^2 ) ^{\varepsilon /2 } ) v^{\prime} \right|
\\
\nonumber
 & \leq & C_1 + \frac{1 }{q+1 } \left| \int_{\mathbb{S} ^{N-1} } v^{q+1 } \Big| _1 ^R \right| + \left|
\int_1 ^R \int_{\mathbb{S} ^{N-1} } O( r^{-1-\varepsilon } + r ^{-1- \gamma \varepsilon } ) \right|
\\
\label{Jest}
 & \leq & C_2 ,
\end{eqnarray}
for some positive constants $C_1 ,C_2  $, independent of $R$. This proves (\ref{estv1}).
\\
To show (\ref{estv2}), we multiply (\ref{sfer1}) with $r^3 v^{\prime \prime} $ and integrate over
$(0,R) \times \mathbb{S} ^{N-1}$ to obtain
\begin{eqnarray}
\nonumber
K & := & \int_0 ^R \int_{\mathbb{S} ^{N-1} } r^{\gamma + 3 } f(r, r^{-\gamma } v) v^{\prime \prime }
\\
\nonumber
 & = & - \int_0 ^R \int_{\mathbb{S} ^{N-1} } r^3 (v^{\prime \prime}) ^2 - (N-1 -2\gamma )
 \int_0 ^R \int_{\mathbb{S} ^{N-1} } r^2 v^{\prime } v^{\prime \prime}
\\
\label{rtheta2}
 & & + \gamma (N-2 -\gamma ) \int_0 ^R \int_{\mathbb{S} ^{N-1} } r v v^{\prime \prime }
-  \int_0 ^R \int_{\mathbb{S} ^{N-1} } r v^{\prime \prime} \Delta _{\theta } v.
\end{eqnarray}
Using integration by parts, we find that the right-hand side of (\ref{rtheta2}) equals to
\begin{eqnarray*}
 & & - \int_0 ^R \int_{\mathbb{S} ^{N-1} } r^3 (v^{\prime \prime }) ^2 - \int_0 ^R \int_{\mathbb{S} ^{N-1} }
 r|\nabla _{\theta } v^{\prime }| ^2
\\
 & & + [N-1- \gamma N +\gamma ^2 ] \int_0 ^R \int_{\mathbb{S} ^{N-1} } r (v^{\prime }) ^2
- (N-1 -2\gamma ) \int_{\mathbb{S} ^{N-1} } \frac{r^2 (v^{\prime}) ^2}{2} \Big| _0 ^R
\\
 & & + \gamma (N-2-\gamma ) \int_{\mathbb{S} ^{N-1} } [ r v v^{\prime } - (1/2) v^2 ] \Big| _0 ^R
\\
 & &
+ \int_{\mathbb{S} ^{N-1} } [ r \nabla _{\theta }
v \cdot \nabla _{\theta } v^{\prime }  -(1/2)|\nabla _{\theta } v|^2  ]  \Big| _0 ^R  .
\end{eqnarray*}
The third term and the last three terms are uniformly bounded in $R$ 
by (\ref{estv1}) and (\ref{ineq3}). Since the second term is nonpositive,
it remains to show that the term $K$ on the left-hand side of (\ref{rtheta2})
is uniformly bounded in $R$. Using again integration by parts, it gives
\begin{eqnarray}
\nonumber
K & =  &  \int_{\mathbb{S} ^{N-1} } r^{\gamma +3} f(r, r ^{-\gamma } v) v^{\prime} \Big| _0 ^R
-(\gamma +3 ) \int_0 ^R \int_{\mathbb{S} ^{N-1} } r^{\gamma +2 } f(r, r ^{-\gamma } v) v^{\prime}
\\
\nonumber
 & & -   \int_0 ^R \int_{\mathbb{S} ^{N-1} } r^{\gamma +3 } f_r (r, r ^{-\gamma } v) v^{\prime}
\\
\nonumber
 & & - \int_0 ^R \int_{\mathbb{S} ^{N-1} }  f_u (r, r ^{-\gamma } v)
 ( -\gamma r^2 v v^{\prime} + r ^3 (v^{\prime} ) ^2 )
\\
\label{estK1}
 & =: & K_1 + K_2 + K_3 + K_4 .
\end{eqnarray}
Using the assumptions on $f$ and (\ref{ineq3}), and proceeding
analogously as in the estimate (\ref{Jest}), we confirm that the terms
$K_i $, ($i=1, 2,3$)
are uniformly bounded in $R$.
Finally, we estimate the term $K_4 $.  By (\ref{fxluq}), we have 
\begin{eqnarray*}
\nonumber
K_4  & = & \int_0 ^R \int_{\mathbb{S} ^{N-1}} q\left( \gamma v^q v^{\prime} -r v^{q-1} (v^{\prime})^2 \right) 
\\
 & & 
+
\int_0 ^R \int_{\mathbb{S} ^{N-1}} q\left( \gamma v^q v^{\prime } -r v^{q-1} (v^{\prime})^2 \right) 
\cdot O (( r^{-2} + r^{-2\gamma } v^2 )^{\varepsilon /2})
\\
 & =: & L_1 + L_2 .
\end{eqnarray*} 
Then, using (\ref{ineq3}), and (\ref{estv1}), it follows that
\begin{eqnarray*}
|L_1 | & \leq  & C_1  \int_{\mathbb{S} ^{N-1} } v^{q+1} \Big| _0 ^R    + C_2 
\int_0 ^R \int_{\mathbb{S} ^{N-1} } r (v^{\prime})^2 \leq C_3 ,
\\
|L_2 | & \leq & C_4 + C_5 \int _1 ^R  \int_{\mathbb{S} ^{N-1} }
r^{-1-\varepsilon } + r^{-1 -\varepsilon \gamma } ) \leq C_6 ,
\end{eqnarray*}
with constants $C_1, \ldots , C_6  $, independent of $R$. Hence $K_4 $ is uniformly bounded in $R$.
This proves (\ref{estv2}).
\\
{\sl Step 2 : } We claim
\begin{equation}
\label{limint}
\lim_{r\to \infty } r^2 \int_{\mathbb{S} ^{N-1 } } (v^{\prime } ) ^2 (r, \theta )\, d\theta =0.
\end{equation}
Suppose by contradiction that (\ref{limint}) is not true.
Then, there exists a sequence $\{ r_n \} $ with
$\lim _{n\to \infty } r_n = +\infty $, and a constant $c>0 $ such that
$$
r_n ^2 \int_{\mathbb{S} ^{N-1 } } (v^{\prime } ) ^2 (r_n , \theta )\, d\theta =: c_n \geq c , \quad n=1,2, \ldots
$$
By (\ref{ineq3}), we have for $r>0$,
\begin{equation*}
\left| \left( r^2  \int_{\mathbb{S} ^{N-1 } } (v^{\prime} ) ^2  \right) ^{\prime} \right|
  =  \left| 2r \int_{\mathbb{S} ^{N-1 } } (v^{\prime} ) ^2 + 2r^2 \int_{\mathbb{S} ^{N-1 } } v^{\prime} v^{\prime \prime }  \right|
 \leq  \frac{M}{r} ,
\end{equation*}
for some $M>0$ independent of $r$.
Hence,
$$
r^2 \int_{\mathbb{S}^{N-1}} (v^{\prime } (r,\theta ))^2 
  \geq  c_n -M \log \frac{r}{r_n }  \geq c
-M \log \frac{r}{r_n } \qquad \forall r>r_n .
$$
In particular,
$$
r^2 \int_{\mathbb{S}^{N-1}} (v^{\prime } (r,\theta ))^2  \geq \frac{c}{2} 
\qquad \forall r\in \left( r_n , r_n e^{c/(2M)} \right) .
$$
But this implies
$$
r \int_{\mathbb{S} ^{N-1 } }
(v' (r, \theta ) ^2 \, d\theta \not\in L^1 (0, +\infty ),
$$
a contradiction. Hence (\ref{limint}) holds.
\\
Furthermore, in view of (\ref{ineq3}), the family $\{ r v^{\prime } (r, \cdot )\} $, $r>0$, is equi--continuous and uniformly in $C^{\tau } (\mathbb{S} ^{N-1}) $
for every $\tau \in (0,1)$.
Denote by ${\mathcal Y} $ the limit set of
$ \{ r v' (r, \cdot )\} $, as $r\to +\infty $.
We claim that ${\mathcal Y} = \{0 \} $. Indeed, let $\omega  \in {\mathcal Y} $.
Then, there exists a sequence $\{ r_n v^{\prime } (r_n , \cdot ) \} $
converging to $\omega $ uniformly in
$C^{\tau } (\mathbb{S} ^{N-1}) $.
By the Dominated Convergence Theorem and (\ref{limint}), we have
$$
\int_{ \mathbb{S} ^{N-1 } }
\omega ^2 (\theta )\, d\theta = \lim _{n\to \infty } r_n ^2
\int_{ \mathbb{S} ^{N-1 } } (v^{\prime } (r_n ,\theta ) ) ^2 \, d\theta =0.
$$
Therefore $\omega (\theta )\equiv 0 $, that is, ${ \mathcal Y } = \{ 0 \} $. In particular,
(\ref{limv1}) holds.
\\
The proof of (\ref{limv2}) is analogous and will be omitted.
$\hfill \Box $
\\[0.1cm]
{\bf Remark 3.8. } The method used in the last proof does not work in the critical case $q= q_2 (\ell)$. This is reminiscent to the fact that problem (\ref{henon1}) possesses solutions of the type (\ref{critsol2}) when $\Omega = \mathbb{R}^N \setminus \{ 0\} $, which do not satisfy (\ref{liminfinity}).
\\[0.1cm]
{\bf Theorem 3.9. } {\sl  Assume that one of the conditions {\bf (i)} or {\bf (ii)} holds. 
Then $u$ satisfies either 
\begin{equation}
\label{limpreciseL}
\lim_{|x|\to \infty } |x|^{\gamma } u(x) =L,
\end{equation}
where 
\begin{equation}
\label{valueofL}
L= \left( \gamma (N-2-\gamma) \right) ^{1/(q-1)} ,
\end{equation}
or 
\begin{equation}
\label{limprecise0}
\lim_{|x|\to \infty } |x|^{\gamma } u(x) =0.
\end{equation}
Moreover,  in the case {\bf (ii)} only  (\ref{limpreciseL}) is possible.
Finally, the convergence is uniform in
$C^{2, \tau } (\mathbb{S} ^{N-1} )$
for any $\tau \in (0, 1)$.}
\\[0.1cm]
To prove Theorem 3.9, we will need a  uniqueness result for non-linear elliptic equations on the sphere. It has been given in \cite{BVV}, Theorem 6.1, (for a slightly weaker result see Corollary B1 and B2 of \cite{GidSpruck}).
\\[0.1cm]
{\bf Lemma 3.10. } {\sl
Let $N\geq 3$, $a >0 $ and $q>1$, and let $V$ be a solution of
\begin{eqnarray}
\nonumber
 & & V\in C^2 (\mathbb{S} ^{N-1} ),
\\
\label{sphericallaplace}
 & & -\Delta _{\theta } V +a V= V^q , \ V>0, \quad \mbox{on }\ \mathbb{S} ^{N-1} ,
\end{eqnarray}
where $\Delta _{\theta } $ denotes the Laplace--Beltrami operator on
$\mathbb{S} ^{N-1}$. Furthermore, assume that 
\begin{equation}
\label{aq1} 
a  \leq  \frac{N-1}{q-1},
\end{equation}
and if $ N\geq 4 $, then also
\begin{equation}
\label{aq2}
q\leq q_S .
\end{equation} 
Finally, assume that one of the inequalities (\ref{aq1}), (\ref{aq2}) is strict in the case $N\geq 4 $.
Then 
\begin{equation}
\label{Vconst}
V\equiv a ^{1/(q-1)} \quad \mbox{on } \ \mathbb{S} ^{N-1} .
\end{equation}
}
{\sl Proof of Theorem 3.9 :}  Let $\{r_n \} $
be a sequence with $\lim _{n\to \infty } r_n = + \infty $, and  let $v$ be defined by (\ref{defv}).
Setting $v_n (\theta ) := v(r_n , \theta ) $, ($\theta \in \mathbb{S} ^{N-1}$),
we have $v_n \in C^{2, \tau } (\mathbb{S} ^{N-1}) $ for every $\tau \in (0,1)$.
By the estimates (\ref{ineq3})
there is a subsequence $\{v_n ' \} $ converging uniformly in
$C^{2,\tau } (\mathbb{S} ^{N-1})$ to a limit $V =V (\theta )$.
Notice
\begin{equation}
\label{estf}
r_n  ^{\gamma +2 } f(r_n , (r_n)^{-\gamma } v_n ) = (v_n ) ^q +
O (r _n ^{-\varepsilon \min \{ 1;  \gamma \} } ) 
\quad \mbox{ as $n\to \infty $.}
\end{equation}
Multiplying (\ref{sfer1}) by $r_n ^2 $, letting $n\to \infty $,
passing to a subsequence,
and taking into account
the assumptions (\ref{fxluq}) and the estimates (\ref{estv1}),
(\ref{estv2}), and (\ref{estf}), we obtain that $V$ satisfies
the equation 
\begin{equation}
\label{Veqn}
-\Delta _{\theta } V + a V =  V ^q \quad \mbox{on }\
\mathbb{S} ^{N-1} ,
\end{equation}
where
\begin{equation}
\label{a}
a:= \gamma (N-2-\gamma ).
\end{equation}
Since $u$ is positive, we have $ V\geq 0 $.
We shall show that $V$ is constant, that is,
\begin{equation}
\label{V=L}
V\equiv a^{1/(q-1)} \quad \mbox{or }\quad V\equiv 0.
\end{equation}
In case {\bf (ii)} we have that $q<1$ and (\ref{V=L}) follows from the uniqueness of the solution of (\ref{Veqn}). Moreover, in view of inequality (\ref{lowerboundii}) of Lemma 3.3, only the first alternative in (\ref{V=L}) is possible. 
\\
Next, consider the case {\bf (i)}.  Then $q>1$, so that we may apply Lemma 3.10. 
\\
Assume first that $N=3 $. Then,
$$
2-a(q-1)  =  2-(2-\ell) \left( 1-\frac{2-\ell}{q-1} \right)
 >  \ell \geq 0,
$$
and (\ref{aq1}) follows with strict inequality. 
\\
Next, assume that $N\geq 4$. We define 
$$
\varphi (z) := \frac{N-3}{4} z^2 - (N-2)z +N-1  \quad z\in (0,2].
$$
It is easy to see that $\varphi (z) >0$ on $(0,2)$ and $\varphi (2)=0$. Hence we have
\begin{eqnarray}
\label{chain}
  N-1 -a (q-1)  &=&  
 N-1 -(2-\ell) (N-2) + (2-\ell)^2 \frac{1}{q-1} 
\\
\nonumber
 & \geq &  N-1 -(2-\ell) (N-2) + (2-\ell)^2 \frac{N-3}{4 } 
\\
\nonumber
 &=&  \varphi (2-\ell)  \geq  0. 
\end{eqnarray}
Moreover, one of the inequalities in the chain (\ref{chain}) is strict in either one of the cases $\ell >0$ and $q\leq q_S $, or $\ell =0$ and $q<q_S $. This means that (\ref{aq1}) follows with strict inequality. 
 Now the Theorem follows from Lemma 3.10.
$\hfill \Box $
\section*{4. A maximum principle for open sets in a half space}
\setcounter{section}{4}
\setcounter{equation}{0}
In this section we obtain a maximum principle for open sets contained 
in a halfspace. 
We will make use of a general comparison principle associated to elliptic operators. A version for bounded domains can be found in \cite{ProtWein}, Theorem 10, and an extension to the case of unbounded domains has been obtained in  
\cite{Naito}, Lemma A. We present the last result in a slightly different form. For the convenience of the reader, we include the full proof.
\\[0.1cm] 
\hspace*{0.3cm}Let $\Omega $ be a domain in $\mathbb{R}^N $ and assume that $a_{ij}$, $b_i $ and $c$ are locally bounded 
measurable functions on $\Omega $ with
\begin{equation}
\label{ellcoeff}
 c_0 |\xi |^2 \leq \sum_{i,j=1}^N a_{ij} (x)\xi _i \xi _j \leq C_0  |\xi |^2 , \qquad 
\xi \in \mathbb{R} ^N ,\quad c_0 ,\ C_0 >0.
\end{equation}
{\bf Lemma 4.1. } {\sl 
Suppose that $u\in W^{2,N} _{loc} (\Omega )\cap C(\overline{\Omega })$ satisfies 
\begin{eqnarray}
\label{uineq1}
 & & -\sum_{i,j=1} ^N a_{ij}  u_{x_i x_j } +\sum_{i=1}^N b_i u_{x_i } 
\leq cu \quad \mbox{ in }\ \Omega ,
\\
\label{uineq2}
 & & u\leq 0 \quad \mbox{on }\ \partial \Omega .
\end{eqnarray}
Suppose, furthermore, that there exists a function 
$w\in W^{2,N} _{loc}(\Omega )\cap C(\overline{\Omega })$ 
such that 
\begin{eqnarray}
\label{wineq1}
 & & w > 0 \quad \mbox{on }\ \overline{\Omega } ,
\\
\label{wineq2}
 & & -\sum_{i,j=1} ^N a_{ij}  w_{x_i x_j } +\sum_{i=1}^N b_i w_{x_i } \geq cw   \quad \mbox{in }\ \Omega .
\end{eqnarray}
Finally, if $\Omega $ is unbounded, we add the requirement that
\begin{eqnarray}
\label{u/wlim}
 & & \limsup _{|x|\to \infty ,\, x\in \Omega } \, 
\frac{u (x )}{w(x )} \leq 0 . 
\end{eqnarray}
Then $u \leq 0 $ in $\Omega $.
}
\\[0.1cm]
{\sl Proof: }  Assume to the contrary that
$u(x^0) > 0$ for some $x^0 \in \Omega $. Choose $\delta  > 0$ such that 
$u(x^0) -\delta w(x^0) = 0$. Define $\overline{u}:= u-\delta w$. In view of (\ref{u/wlim}), there exists $R > |x^0|$ such that $ \overline{u}\leq  0$ on $\partial B_R \cap \Omega $. Then 
$$ 
-\sum_{i,j=1} ^N a_{ij}  \overline{u}_{x_i x_j } +\sum_{i=1}^N b_i \overline{u}_{x_i } \leq c\overline{u}
\quad \mbox{
on $\Omega \cap  B_R$ }
$$ 
and $\overline{u}\leq 0 $  on $\partial ( \Omega \cap  B_R)$. 
Setting $v:= \overline{u}/w $, a short calculation then shows that 
\begin{eqnarray*}
 & & -\sum_{i,j=1} ^N a_{ij}  v_{x_i x_j } +\sum_{i=1}^N b_i v_{x_i } - \frac{2}{w} \sum_{i,j=1} ^N   a_{ij} w_{x_j}  v_{x_i } \leq 0
\quad \mbox{
on $\Omega \cap  B_R$, } 
\\
 & & v\leq 0 \quad \mbox{on }\ \partial (\Omega  \cap B_R ),
\end{eqnarray*}
while $v(x^0 )=0$. 
By the Strong Maximum Principle (Theorem 7.1, {\bf (I)}), this implies that $v\equiv 0 $ in $\Omega \cap B_R $. 
But this means that $u= \delta w >0 $ on $\partial (\Omega \cap B_R )$, a contradiction.
$\hfill \Box $
\\[0.1cm] 
\hspace*{0.3cm}Appropriate choices of the comparison function $w$ in Lemma 4.1 lead to the following maximum principle for open subsets of a halfspace.  It will enable us to treat cases of slow decay in the proof of radial symmetry.
\\[0.1cm]
{\bf Theorem 4.2. } {\sl
Let $\Omega $ be an open set with $\overline{\Omega } \subset \mathbb{R}^N_+$, $a\in [0, N/2) $, $b\geq 0$, and let $u\in W^{2,N} _{loc} (\Omega ) \cap C (\overline{\Omega } )$ be such that
\begin{eqnarray}
\label{bound1}
 & &  u(x) \leq C |x|^{-b} , 
\\
\label{diffineq1}
 & & -\Delta u \leq K \frac{u}{|x|^2 }  \quad \mbox{ in $\Omega $,}
\\
\label{bdry}
 & & u\leq 0 \quad \mbox{on } \ \partial \Omega ,
\end{eqnarray}
for some positive constants $C$ and $K$, where
\begin{eqnarray}
\label{bound2}
 & & K \leq  \frac{N^2 }{4} -a^2 \ \ \mbox{ and}
\\
\label{bound3}
 & & N<2a +2b+2.
\end{eqnarray}   
Then $u\leq 0 $ in $\Omega $.
}
\\[0.1cm]
{\sl Proof : } Defining
\begin{equation}
\label{wdefine}
w(x) := x_1 |x|^{a-(N/2)} , \quad x\in \Omega,
\end{equation}
we have
\begin{equation}
\label{wequation}
 -\Delta w = \left( \frac{N^2 }{4} -a^2 \right) \frac{w}{|x|^2}  \quad w>0 \quad \mbox{in }\ \overline{\Omega}. 
\end{equation} 
Next, for $R>0$, we set $\Omega _R := \{ x:\ Rx \in \Omega \} $ and
$$
u_R (x) := R^b u(Rx) , \quad x\in \Omega _R. 
$$
Since $\Omega _R \subset \mathbb{R}^N _+ $,
\begin{eqnarray*} 
 & & 0\leq u_R (x) \leq C|x| ^{-b } , \quad -\Delta u_R (x)\leq C K|x|^{-b-2} \ \mbox{ in $\ \Omega _R \ $  and}
\\
 & & u_R (x)\leq 0 \ \mbox{ on $\ \partial \Omega _R $ }.
\end{eqnarray*}
Standard elliptic estimates show that $u_R (x) \leq D x_1 |x| ^{-b-1} $ in $\Omega _R \cap (B_2 \setminus \overline{B_1 } )$, where $D$ is a positive constant that does not depend on $R$. This also implies that  $u (x) \leq D x_1 |x| ^{-b-1} $ in $\Omega  \cap (B_{2R} \setminus \overline{B_R } )$.
Since $R$ was arbitrary, it follows that
\begin{equation}
\label{improvedboundu}
u(x)\leq Dx_1 |x|^{-b-1} \quad \mbox{in $\Omega $.}
\end{equation}  
Now 
the assertion follows from Lemma 4.1, taking  
 $a_{ij} =\delta _{ij} $, $b_i =0$, ($i,j=1, \ldots ,N$), and 
$c(x) = K$.
$\hfill \Box $   
\section*{5. Moving Plane Method and symmetry }
\setcounter{section}{5}
\setcounter{equation}{0}
In this section we use the
Alexandroff--Serrin Moving Plane Method 
to show the Theorems  1.1 and 1.2.  
We will need the following two technical Lemmata.
\\[0.1cm]
{\bf Lemma 5.1. } 
{\sl 
Assume that one of the conditions {\bf (i)} or {\bf (ii)} 
holds. 
Then there exists a number $a\in [0, N/2) $ such that
\begin{eqnarray}
\label{a1}
 & &  N<2a +2\gamma +2 \quad \mbox{and }
\\
\label{a2}
 & &  \max \{ 1,\, q\} \cdot \gamma (N-2-\gamma ) <\frac{N^2 }{4} -a^2 .
\end{eqnarray}
}
{\sl Proof :} 
\\
{\sl Case {\bf (i)} }: Assume first that $q\in (q_1 (\ell)  , q_2 (\ell) )$. Then 
\begin{eqnarray*}
N & < & 2\gamma +2 \quad \mbox{and }
\\
q \gamma (N-2-\gamma ) & < & q_2 (\ell)  \left( \frac{N-2}{2} \right) ^2 
= \frac{(N+2-2\ell ) (N-2)}{4} 
 <  \frac{N^2 }{4}.
\end{eqnarray*}
Hence (\ref{a1}) and (\ref{a2}) are satisfied with $a=0$.  
\\
Now let $q > q_2 (\ell)$. We set
\begin{eqnarray}
\label{defa0}
a_0 & := & \frac{N-2}{2} -\gamma \quad \mbox{ and}
\\
\label{defD}
D & := & \frac{N^2 }{4}-a_0 ^2 -q \gamma (N-2-\gamma ) 
\\
\nonumber 
 &= & (2-\ell) \left(
\frac{2-\ell}{q-1} -N+2\right) + N-1.
\end{eqnarray}
 Note that $a_0 \in (0, N/2)$ and (\ref{a1}) is satisfied for every $a\in (a_0 , N/2] $. 
\\
Now, if $N=3 $, we have  
$$
D  =  \frac{(2-\ell )^2 }{q-1 } + \ell >0. 
$$
Further, if  $N\geq 4$, then we obtain,  
\begin{eqnarray}
\label{chain1}
D & \geq &  (2-\ell) \left(
\frac{2-\ell}{q_S  -1} -N+2\right) + N-1  
\\
\nonumber 
 & = &    
 (2-\ell) \left( \frac{(2-\ell)(N-3)}{4} -N+2 \right) + N-1 
\\
\nonumber 
 & \geq & 0 .
\end{eqnarray} 
Moreover, one of the inequalities in the chain (\ref{chain1}) is strict in either one of the cases $\ell >0$ and $q\leq q_S $, or $\ell =0$ and $q<q_S $.  
\\
Since $D>0$ in both cases, 
 we find $a\in (a_0 , N/2]$ with $|a-a_0 | $ small, such that both (\ref{a1}) and (\ref{a2}) hold.
\\
{\sl Case  {\bf (ii)}} : Here we have $q <q_1 (\ell)  <1 $ and $\max \{ 1,\, q \} =1 $. 
First, assume that 
$$
\max \{ 0,\,  q_2 (\ell) \} <q<  q_1 (\ell)  .
$$ 
Then 
$$
 \gamma (N-2-\gamma ) < \left(\frac{N-2}{2}\right) ^2 
 <\frac{N^2 }{4} \quad \mbox{and }
\ \
\gamma  =  \frac{\ell-2}{1-q} > \frac{\ell-2}{1-q_2 (\ell)  } =  \frac{N-2}{2},
$$
then (\ref{a1}) and (\ref{a2}) are satisfied with $a=0$.
\\
Next assume that $q_2 (\ell)>0$ and $q\in (0,q_2 (\ell)  ]$. Let $a_0 $ again be given by 
(\ref{defa0}).
Then we have $a_0 \in [0, N/2) $, 
and (\ref{a1}) is satisfied for every $a\in (a_0 ,N/2 )$. 
Notice that there holds 
$$
\frac{N^2 }{4}-a_0 ^2 -
\gamma (N-2-\gamma )  =   
 N-1 >0.
$$
Hence we can find $a\in (a_0 , N/2) $ with $|a-a_0 |$ small, such that both (\ref{a1}) and (\ref{a2}) hold.
$\hfill \Box$
\\[0.1cm]
{\bf Lemma 5.2. } {\sl Assume that (\ref{fxluq}) and one of the conditions {\bf (i)}  or {\bf (ii)} are satisfied. Furthermore, let $u$ be a solution of {\bf (P)} or {\bf (P)}$_{\mathbf{0}}$.  Then, in the case {\bf (i)}  there holds either
\begin{equation}
\label{basiclimit1}
\lim_{|x|\to \infty } |x|^{2-\ell} u^{q-1} (x) = \gamma (N-2-\gamma ),
\end{equation} 
or 
\begin{equation}
\label{basiclimit2}
\lim_{|x|\to \infty } |x|^{2-\ell} u^{q-1}  (x) = 0.
\end{equation}  
In the case {\bf (ii)} there holds only (\ref{basiclimit1}). }
\\[0.1cm]
{\sl Proof: } Using Theorem 3.9 we have that
$$
\lim_{|x|\to \infty } |x|^{2-\ell} u^{q-1} (x) = \left( \lim_{|x|\to \infty } |x|^{\gamma } u(x)\right)^{q-1} = L^{q-1},
$$
and the assertions follow.
$\hfill \Box $
\\[0.1cm]
{\bf Lemma 5.3. } {\sl Assume that (\ref{fxluqsimple}) and one of the conditions {\bf (iii)}--{\bf (vi)} are satisfied.
Furthermore, let $u$ be a solution of problem {\bf (P)} or of {\bf (P)}$_{\mathbf{0}}$.
Then the limit property (\ref{basiclimit2}) holds.}
\\[0.1cm]
{\sl Proof: } First observe that (\ref{basiclimit2}) is trivial in the case {\bf (iii)}.
Furthermore, if one of the conditions {\bf (v)} or {\bf (vi)} is satisfied, then we have by the Lemmata  3.1 and 3.4
$$
c_1 |x|^{2-N} \leq u(x) \leq c_3 |x|^{2-N} \quad \mbox{on $B_1 ^C $},
$$
which implies that 
$$
|x|^{2-\ell} u^{q-1}  \leq  
C |x|^{N-\ell -q (N-2)} 
$$
for these $x$, for some positive constant $C$. Since $N-\ell -q(N-2)<0$, we deduce (\ref{basiclimit2}).
\\
Finally, if {\bf (iv)} is satisfied, then $q= q_1 (\ell) <1$ and $u$ satisfies (\ref{frombelow}), so that with some positive constant $C'$,
\begin{eqnarray*}
|x|^{2-\ell} u^{q-1} & \leq & C' |x|^{2-\ell} \left[ |x|^{2-N} \log |x| \right] ^{q-1}
\\
 & = & C' \left[ \log |x|\right] ^{q-1}  \qquad \forall x\in \mathbb{R}^N \setminus{B_2 } ,
\end{eqnarray*}
and (\ref{basiclimit2}) follows again. 
$\hfill \Box$
\\[0.1cm] 
\hspace*{1cm}Now we are in a position to prove our symmetry results.  Let us introduce some classical notation. 
For $\lambda \in \mathbb{R}$, let
\begin{eqnarray*}
T^{ \lambda }   & = & \{ x= (x_1 , \ldots , x_N ) \in \mathbb{R}^N :\, x_1 =\lambda \} ,
\\
\Sigma (\lambda )   & = & \{ x= (x_1 , \ldots , x_N )\in \mathbb{R}^N :\, x_1 >\lambda \} .
\end{eqnarray*}
For $x\in \mathbb{R}^N $, let $x^{\lambda } $
denote the reflection point of $x$ about $T^{\lambda } $, that is,
$$
x^{\lambda } := (2 \lambda -x_1 , x_2 ,\ldots , x_N ).
$$
Further, 
let
\begin{eqnarray*}
u^{\lambda } (x) & := &  u(x^{\lambda } ), \quad (x\in \mathbb{R}^N ),
\\
w^{\lambda } & := & u-u^{\lambda } ,
\\
\Omega (\lambda ) & := & 
\{x \in \Sigma (\lambda ):\, w^{\lambda } (x) >0 \} , \qquad \mbox{and}
\\
A & := & \{ \lambda > 0 :\, w^{\lambda } (x) <0   \ \ \forall x\in \Sigma (\lambda ) \} .
\end{eqnarray*}
{\sl Proof of Theorem 1.1 for Problem {\bf (P)}:}  
\\
We proceed in $7$ steps.
\\
{\sl Step 1: }
Define
$$
c^{\lambda } (x):= 
\left\{
\begin{array}{ll}
\frac{ f(|x|, u(x) ) -f(|x|, u^{\lambda } (x)}{u(x)-u^{\lambda } (x) } & 
\mbox{ if }\
u(x)\not= u^{\lambda } (x), 
\\
0 & \mbox{ if } 
\ u(x)= u^{\lambda } (x).
\end{array}
\right.
$$
By the assumptions on $f$ the functions $c^{\lambda } $ are locally bounded in $\Sigma (\lambda )$. Furthermore, if $\lambda \geq 0$, then we have $|x|\geq |x^{\lambda }|$ for all $x\in \Sigma (\lambda)$, so that condition (\ref{fradial}) gives  
\begin{eqnarray}
\nonumber
 - \Delta w^{\lambda}
 & = &  
f(|x|, u) - f(|x^{\lambda } | , u^{\lambda } )
\\
\label{wlambda1}
 & \leq & 
f(|x|, u) - f(|x| , u^{\lambda } ) 
= c^{\lambda } w^{\lambda } \qquad \mbox{on} \ \Sigma (\lambda ) \ {\rm if} \ \lambda \geq 0.
\end{eqnarray}
Furthermore, we can find $R_0 >0$, such that 
\begin{equation}
\label{f>0}
 f(|x|, u(x)) >0 \quad \mbox{for $|x|\geq R_0 $.}
\end{equation} 
Since $u$ is positive and $\lim_{|x|\to \infty } u(x)=0$, we may add the following requirements for $R_0 $:
\begin{eqnarray}
\label{maxuleft}
 & & \max \{ u(x): \, x_1 \leq R_0 \} > \max \{ u(x):\, x_1 \geq R_0 \}  \ \mbox{ and }
\\ 
\label{maxuright}
 & & \max \{ u(x): \, x_1 \geq - R_0 \} > \max \{ u(x):\, x_1 \leq -R_0 \} .
\end{eqnarray} 
The last property (\ref{maxuright}) implies that 
\begin{equation}
\label{Abddbelow}
A\subset (-R_0 , + \infty).
\end{equation}
Moreover, since $w^{\lambda } < 0 $ and satisfies (\ref{wlambda1}) in $\Sigma (\lambda ) $ if $\lambda \geq 0 $ and $\lambda \in A$, Hopf's Boundary Point Lemma (Theorem 7.1, {\bf (II)}) yields  
\begin{equation}
\label{hopf} 
\frac{\partial}{\partial x_1 } w^{\lambda } (x) <0 \quad \mbox{on  $T^{\lambda } $ \ if $\lambda \in A\cap [0,+\infty ) $.} 
\end{equation}
{\sl Step 2 : } 
Next, we estimate the functions $c^{\lambda } (x)$ on the sets $\Omega (\lambda )\setminus \overline{B_{R_0 } }$ if $\lambda \geq 0$. 
We claim that there is a function
$m:[0, +\infty ) \to [0, + \infty )$, with $\lim _{s\to 0} m(s)=0$, such that
\begin{equation}
\label{mainineq}
c^{\lambda } \leq \max\{ 1,\, q\} \cdot |x|^{-\ell} u^{q-1} 
\left(
1+ m(u^2 + |x|^{-2}) 
\right) 
\quad 
\forall x\in \Omega (\lambda )\setminus \overline{B_{R_0 } }.
\end{equation} 
To prove (\ref{mainineq}) we split into two cases.
\\
{\sl Case {\bf (i)}: } Then we have   $q>1$. Using the assumptions (\ref{fradial}) and (\ref{fxluq}) on $f$ and defining 
$$
h_t := u	^{\lambda } +t w^{\lambda }  \quad (t\in [0,1]),
$$
we obtain
\begin{eqnarray}
\nonumber 
f(|x|, u)-f(|x|, u^{\lambda } ) 
 & = &  
w^{\lambda } \int _0 ^1 f_u (|x|, h_t )\, dt 
\\
\nonumber
 &  \leq & 
q|x|^{-\ell} w^{\lambda } \int _0 ^1 
(h_t)^{q-1} 
\left( 
1+ m_1 ((h_t)^2 + |x|^{-2} )
\right) 
\, dt 
\\
\label{aux1}
\end{eqnarray}
$\forall x\in \Omega (\lambda ) \setminus \overline{B_{R_0 }}$, $m_1:[0,+\infty) \to  [0,+\infty) $ is an increasing  function with 
$\lim _{s\to 0} m_1 (s) =0$. 
Furthermore, since $w^{\lambda } >0$ in 
$ \Omega (\lambda)$, we have that
$h_t \leq u$. Thus, (\ref{aux1}) implies that
\begin{eqnarray}
\nonumber 
 f(|x|, u)-f(|x|, u^{\lambda } )
  & \leq &  
q|x|^{-\ell} u^{q-1} 
\left( 
1+ m_1 (u^2 + |x|^{-2} )
\right) 
w^{\lambda } 
\\
\label{ineqcase1} 
 & &
\quad \forall x\in \Omega (\lambda ) \setminus \overline{B_{R_0 } }. 
\end{eqnarray}
{\sl Case {\bf (ii)}: } 
Then we have $q<1$. With the notation from the previous case, we obtain on $\Omega (\lambda )\setminus \overline{B_{R_0 } }$,
\begin{eqnarray}
\nonumber
f(|x|, u) - f(|x|, u^{\lambda })
 & \leq & 
u^{q-1} \left( 
f(|x|, u) u^{1-q} - f(|x|, u^{\lambda } ) (u^{\lambda } )^{1-q} 
\right)
\\
\label{q<1}
 & = & u^{q-1} w^{\lambda } \int_0 ^1 \frac{d}{du} \left[ f(|x|, h_t  ) (h_t) ^{1-q} \right] \, dt .
\end{eqnarray}
By the assumptions on $f$ and arguing as before, we find an increasing function \\
$m_2 : [0, +\infty ) \to [0,+\infty ) $ with $\lim_{s\to 0} m_2 (s)=0$ such that
\begin{eqnarray}
f(|x|, u)- f(|x|, u^{\lambda }) 
\nonumber
 & \leq &
|x|^{-\ell} u^{q-1} w^{\lambda } \int_0 ^1 \left( 1+ m_2 ((h_t)^2 + |x|^{-2} )  \right) \, dt 
\\
\label{ineqcase2}
 & \leq &
|x|^{-\ell} u^{q-1}  \left( 1+ m_2 (u^2 + |x|^{-2} )  \right) w^{\lambda } 
\quad \mbox{ in $\Omega (\lambda )\setminus \overline{B_{R_0 } } $.}
\end{eqnarray}
Now (\ref{mainineq}) follows from (\ref{ineqcase1}) and (\ref{ineqcase2}). 
Together with  (\ref{wlambda1}), we obtain
\begin{equation}
\label{mainineq2}
-\Delta w^{\lambda } \leq \max \{ 1,\, q\}  \cdot |x|^{-\ell} u^{q-1}  \left( 1+ m (u^2 + |x|^{-2} )  \right) w^{\lambda } \quad \mbox{ on $\Omega (\lambda )\setminus \overline{B_{R_0 } } $.}
\end{equation} 
{\sl Step 3 : } Next, we  apply Theorem 4.2 to the differential inequality (\ref{mainineq2}).  By Lemma 5.1, there is a number $a\in [0,N/2 )$ such that 
(\ref{a1}) holds. We choose $\varepsilon _0 >0$ small enough such that 
\begin{equation}
\label{finer!}
\max \{1,\, q\} \cdot \gamma (N-2-\gamma ) (1+\varepsilon _0 ) < \frac{N^2 }{4} -a^2 .
\end{equation}
Since $u$ satisfies (\ref{decaycrucial})  and either  (\ref{basiclimit1}) or (\ref{basiclimit2}),   
 we may add the requirement to $R_0 $ that 
\begin{equation}
\label{add}
|x|^{-\ell } u^{q-1} \big(1+ m(u^2 + |x|^{-2} ) \Big)
 \leq  \gamma (N-2-\gamma ) ( 1+\varepsilon _0 ) |x|^{-2}  
\quad \forall x \in \mathbb{R}^{N} \setminus \overline{B_{R_0 }} .
\end{equation}
 
Now, using (\ref{mainineq2})-(\ref{add}), we find that 
\begin{equation}
\label{diffineq}
-\Delta w^{\lambda }
\leq 
\left( \frac{N^2 }{4} -a^2 \right) |x|^{-2} w^{\lambda }
\quad \mbox{ in $\Omega (\lambda ) \setminus \overline{B_{R_0 } } $,}
\end{equation}
where
$N<2a +2 \gamma +2$.  
\\
{\sl Step 4 :} We claim:
\begin{equation}
\label{Abegin}
[R_0 , +\infty )\subset A.
\end{equation}
By (\ref{decaycrucial}) we have
\begin{equation}
\label{west1}
w^{\lambda } (x) \leq C |x| ^{-\gamma } \quad \mbox{in }\ \Omega (\lambda ), 
\end{equation} 
for some $C>0$, for every $\lambda \geq  0$.
Furthermore, note that (\ref{diffineq}) holds on the set $\Omega (\lambda )$ for every $\lambda \geq R_0 $. Since $w^{\lambda } = 0$ on $\partial \Omega (\lambda )$ and taking into account (\ref{west1}), Theorem 4.2 tells us that 
$ w^{\lambda } \leq 0 $ on $\Omega (\lambda )$,
which implies that $\Omega (\lambda )=\emptyset $ whenever $\lambda \geq R_0 $. Hence, we have that $w^{\lambda } \leq 0 $ on $\Sigma (\lambda )$ for $\lambda \geq 0$. Now,  assume that there is a $\lambda _0 \in [R_0 , +\infty ) $ and a point $x^0 \in \Sigma (\lambda _0 ) $ with $w^{\lambda _0 } (x^0 )=0$. Then, the Strong Maximum Principle (Theorem 7.1, {\bf (I)}) yields $w^{\lambda _0 } \equiv 0$  in $\Sigma (\lambda _0 )$. But this contradicts (\ref{maxuleft}). Hence, we must have $w^{\lambda } <0$ on $\Sigma (\lambda )$ whenever $\lambda \geq R_0 $. This is (\ref{Abegin}). 
\\
{\sl Step 5 : } We define 
\begin{equation}
\label{lambdastar}
\lambda _+ := \inf \Big\{ \lambda :\, \mu \in A \ \ \forall \mu \in [\lambda , +\infty ) \Big\} .
\end{equation} 
By (\ref{Abegin}) and (\ref{maxuright}), we must have 
\begin{equation}
\label{lambdastarestimate}
\lambda _+ \in [-R_0 , R_0 ] ,
\end{equation}
which also implies that  
\begin{equation}
\label{rightoflambdastar}
 (\lambda _+ , +\infty ) \subset A, 
\end{equation}
and by continuity,
\begin{equation}
\label{wlambdastar}
w^{\lambda _+} \leq 0 \quad \mbox{on }\ \Sigma (\lambda _+ ).
\end{equation}
Since property (\ref{hopf}) holds for all $\lambda >\lambda _+$, it follows that
\begin{equation}
\label{sym1} 
 u_{x_1 }( x) <0 \quad \mbox{on $\Sigma (\lambda _+ ) $.}
\end{equation}
{\sl Step 6 :} We distinguish two cases.
\\
{\bf (a)} Assume that $\lambda_+ > 0$.  
We claim that this implies
\begin{equation}
\label{w=0}
w^{\lambda _+} \equiv 0 \quad \mbox{on $\Sigma (\lambda _+ )$.}
\end{equation}
Suppose that this is not true. Using the Strong Maximum Principle as  in Step 4, we deduce that $w^{\lambda _+} <0$ on $\Sigma (\lambda _+)$. Then, arguing as in Step 1,  Hopf's Boundary Point Lemma tells us that  
\begin{equation}
\label{hopf1}
\frac{\partial }{\partial x_1 } w^{\lambda _+ } (x) <0 \quad 
\mbox{ on $T^{\lambda _+ } $. }
\end{equation}   
By continuity, this implies that there exists a number $\delta >0$ such that $w^{\mu } (x) < 0 $ on 
$\Sigma (\mu ) \cap B_{2R_0 } $, whenever  $0\leq \lambda _+ -\delta \leq  \mu <\lambda _+ $. 
By (\ref{diffineq}), we have
$$
-\Delta w^{\mu }
  \leq  
\left( 
\frac{N^2 }{4} -a^2 
\right) 
|x|^{-2} w^{\mu } \quad 
\mbox{ in $\Omega (\mu ) \setminus \overline{B_{R_0 }} $,}
$$ 
where $N<2a+ 2\gamma + 2$, and (\ref{west1}) with $\mu $ in place of $\lambda $.
Since 
$$
 w ^{\mu }  \leq  0 \quad 
\mbox{ on $\partial \left( \Omega (\mu ) \setminus   \overline{B_{R_0 }} \right) $,}
$$
Theorem 3.2 implies that  $w^{\mu } \leq 0 $ on the set $
\Omega (\mu ) \setminus \overline{B_{R_0}} $.  Using the Strong Maximum Principle, this implies that $w^{\mu } <0 $ in $\Sigma (\mu )$ for these $\mu $. But this contradicts to the definition of $\lambda _+$. Hence (\ref{w=0}) follows, which implies
\begin{equation}
\label{sym2} 
 u( x) = u(x^{\lambda _+ } )  \quad \mbox{on $\Sigma (\lambda _+ ) $.}
\end{equation}
{\bf (b)} Next assume that $\lambda _+ \leq 0$. Defining 
$v(x) := u( -x_1 ,x') $ for $x= (x_1 , x')\in \mathbb{R}^N $, we have that $-\Delta v= f(|x|,v) $. We may then repeat all the above arguments for $v$ in place of $u$. This leads to the existence of a number $\lambda_- $, such that 
\begin{eqnarray}
\label{monv}
 & & u_{x_1 } (x_1 , x' ) >0 \quad \mbox{and }
\\
\label{limv}
 & & u(x_1 , x' ) \leq u(2 \lambda _- -x_1 ,x') \quad \mbox{for $x_1 <\lambda _- $ and $x'\in \mathbb{R}^{N-1} $.}
\end{eqnarray} 
In view of the properties of $u$ that we have already proved, we must have $\lambda_- \leq  \lambda _+$, and in particular, $\lambda_- \leq 0$. Now, if $\lambda_- <0 $, then we conclude as before that  
\begin{equation}
\label{sym2}
u(x_1 , x' ) = u(2\lambda_- -x_1 ,x') \quad \mbox{for $(x_1, x' ) \in  \mathbb{R}^{N} $.}
\end{equation}
On the other hand, if $\lambda_- =0$, then we must also have $\lambda _+ =0 $ and
\begin{equation}
\label{sym3}
u(x_1 , x' ) = u( -x_1 ,x') \quad \mbox{for $(x_1 , x') \in \mathbb{R}^N  $.}
\end{equation}
To sum up, we have proved that there is a number 
$\lambda ^* \in \mathbb{R} $ such that 
\begin{eqnarray}
\label{sym4}
 & & u_{x_1 } (x_1 , x') >0 \quad \mbox{and }
\\
\label{sym5}
 & & u(x_1 , x')= u(2\lambda ^* -x_1 ,x' ) \quad \mbox{for $x_1 <\lambda ^* $ and $x'\in \mathbb{R}^{N-1} $.}
\end{eqnarray}
Since properties (\ref{sym4}) and (\ref{sym5}) hold in every cartesian coordinate system centered at the origin, it follows that
$u$ is radially symmetric and radially decreasing w.r.t. some point $x^0$.
\\ 
{\sl Step 7 : } It remains to prove that $x^0 =0$ if $\ell >0$.  
\\
Assume that this is not the case. Then, there is a coordinate system such that (\ref{sym4}) and (\ref{sym5}) hold with some number $\lambda ^*>0$.   
Putting
$$
\xi (t):= |x^{\lambda ^* }| +t ( |x|- |x^{\lambda ^* }| )  \quad t\in [0,1],
$$ 
and using (\ref{fxluq}), we find  on the set 
$\Sigma (\lambda ^* )$, 
\begin{eqnarray*}
0 &\equiv & -\Delta w^{\lambda ^* } = f(|x|, u)- f(|x^{\lambda ^* }|, u ) 
\\
 & = & (|x|-|x^{\lambda ^* } |)  \cdot \int_0 ^1 f_r (\xi (t), u)\, dt 
\\
 & = & -\ell (   |x|-|x^{\lambda ^* } |) u^q \cdot \int_0 ^1 |\xi (t)| ^{-\ell-1} \left[ 1 + O ((u^2 + |\xi (t)|^{-2} )^{\varepsilon /2}) ) \right] \, dt .
\end{eqnarray*}
Since $|x|>|x^{\lambda ^* }|$, this implies 
$$
 0 =  - \int_0 ^1 |\xi (t)| ^{-\ell-1} \left[ 1 + O ((u^2 + |
\xi (t)|^{-2} )^{\varepsilon/2} ) \right] \, dt .
$$
But this is impossible when $|x| $ is large enough. Hence, we must have $w^{\lambda ^* }<0 $ on $\Sigma (\lambda ^* )$, a contradiction. 
$\hfill \Box$
\\[0.1cm]
{\sl Proof of Theorem 1.4 for Problem {\bf (P)}:}
\\
The proof is analogous to the previous one, except with few modifications in the application of Theorem 4.2 that we detail below.
\\
First, we obtain (\ref{wlambda1}) as before.
Furthermore, we may choose $R_0 >r_0 $ large enough such that $u(x) < u_0 $ for $|x|>R_0 $. Then, applying the assumptions  (\ref{fxluqsimple}), we obtain
\begin{eqnarray}
\label{12a}
 & & 0< f (|x|,u(x))\leq d_1 |x|^{-\ell} (u(x))^q ,
\\
\label{12b}
 & & f_u (|x|, u(x)) \leq d_2 |x|^{-\ell} (u(x))^{q-1} \quad \mbox{for $|x|> R_0 $.} 
\end{eqnarray}
Now we again split into two cases. First assume that $q\geq 1$. 
Since we have
$$ 
u(x)\geq u^{\lambda } (x) +t w^{\lambda } (x)\equiv h_t (x)  \quad \forall
x\in \Omega (\lambda ),
$$ 
it follows that 
\begin{eqnarray}
\nonumber
f(|x|, u) -f(|x|, u^{\lambda }) 
 & = &  w^{\lambda } \int_0 ^1 f_u (|x|, h_t ) \, dt 
\\
\nonumber  & \leq & d_2 |x|^{-\ell} w^{\lambda } \int_0 ^1 (h_t )^{q-1}  \, dt 
\\
\label{12c}
 & \leq & d_2 |x|^{-\ell} u ^{q-1} \quad \mbox{in $ \Omega (\lambda ) \setminus \overline{B_{R_0 }} .$}
\end{eqnarray}
Now let $q< 1$. Then, we again obtain (\ref{q<1}) and the assumptions (\ref{fxluqsimple}) yield 
\begin{eqnarray}
\nonumber
f(|x|, u) -f(|x|, u^{\lambda }) 
 & \leq & w^{\lambda } u^{q-1} \int_0 ^1 \frac{d}{du} \left[ f(|x|,h_t) (h_t)^{1-q} \right] \, dt 
\\
\label{12d}
 & \leq & \left( d_2 + (1-q) d_1 \right) |x|^{-\ell}  u^{q-1 } w^{\lambda } \quad \mbox{in $ \Omega (\lambda ) \setminus \overline{B_{R_0 }} .$}
\end{eqnarray}
Now (\ref{wlambda1}), (\ref{12c}), and (\ref{12d}) show that there is  $d_3>0$, independent of $\lambda $, such that
\begin{equation}
\label{12e}
-\Delta w^{\lambda } \leq d_3 |x|^{-\ell } u^{q-1} w^{\lambda } \quad  \mbox{in $ \Omega (\lambda ) \setminus \overline{B_{R_0 }} .$}
\end{equation}
On the other hand, Lemma 5.3 tells us that 
\begin{equation}
\label{basiclimit3}
\lim_{|x|\to \infty } |x|^{2-\ell} u^{q-1} (x)=0 .
\end{equation} 
Hence, by choosing $R_0 $ large enough in (\ref{12e}), we have that 
\begin{equation}
\label{diffineq2}
-\Delta w^{\lambda } \leq |x|^{-2} w^{\lambda } \quad \mbox{on } \ \Omega (\lambda ) \setminus \overline{B_{R_0} }.
\end{equation}
Furthermore, since $u$ decays at infinity, we also have that
\begin{equation}
\label{decaygeneral}
\lim_{|x|\to \infty ,\ x\in \Omega (\lambda )} w^{\lambda } (x) =0.  
\end{equation}
Now, applying Theorem 4.2, with $a=(N-1)/2 $ and  
$b= 0$,  shows that $w^{\lambda } \leq 0 $ on $\Omega (\lambda ) $ whenever $\lambda \geq R_0 $, which in turn implies  (\ref{Abegin}).  Then, repeating the steps of the last proof and using (\ref{diffineq2}) and (\ref{decaygeneral}) in place of (\ref{diffineq}) and (\ref{west1}), respectively, one proves the symmetry properties (\ref{sym4}) and (\ref{sym5}).  Hence, $u$ is radially symmetric and radially decreasing w.r.t. some point $x^0 \in \mathbb{R}^N $.
$\hfill \Box $
\\[0.1cm]
{\sl Proof of Theorem 1.3 for Problem {\bf (P)}: }
\\
In view of the properties {\bf (i')} and (\ref{faster}), we obtain 
(\ref{basiclimit2}). We may then proceed analogously as in the proof of Theorem 1.3.
$\hfill \Box $ 
\\[0.1cm]
{\sl Proof of the Theorems 1.1, 1.3 and 1.4 for Problem {\bf (P)}$_{\mathbf{0}}$: }
\\
Define 
$z(\lambda ):= (2\lambda , 0,\ldots ,0)$ for $\lambda >0$. Henceforth, we use the notations of the proof of Theorem 1.1, except that we replace in the definition of $A$  the sets $\Sigma (\lambda )$ by
$
\Sigma ' (\lambda ):= \Sigma (\lambda ) \setminus \{ z(\lambda )\} $.   
Observe that we must have $A\subset [0,+\infty )$ since $\lim_{x\to 0} u(x)=+\infty $. If $\lambda >0$, 
we find a number $\varepsilon (\lambda )  \in (0, \lambda ) $ such that $w^{\lambda } (x) <0 $ in $B_{\varepsilon (\lambda )} (z(\lambda ))\setminus \{ z(\lambda )\} $, which means that $B_{\varepsilon (\lambda )} (z(\lambda )) \cap \Omega (\lambda) = \emptyset $. Hence, $w^{\lambda } $ is regular in $\Omega (\lambda )$ and satisfies the differential inequality (\ref{diffineq}). 
Then, defining $\lambda _+ $ by (\ref{lambdastar}) we must have  that $\lambda _+ \geq 0$, and proceeding as before, we obtain that $w^{\lambda _+ } \leq 0$ on $\Sigma ' (\lambda _+ )$ and $u_{x_1 } <0 $ on $\Sigma  (\lambda _+ )$. Now assume that $\lambda _+ >0$. Then, since $w^{\lambda _+ } <0 $ in 
$B_{\varepsilon (\lambda _+ )} (z(\lambda _+ )) \setminus \{ z(\lambda _+ )\}$, the Strong Maximum Principle yields $w^{\lambda _+ } <0$ in $\Sigma '(\lambda _+ )$. But this leads again to a contradiction. Hence, we must have $\lambda _+ =0$. Then, repeating the same analysis for the function $v(x_1 , x'):= u(-x_1 ,x' )$, we find that $w ^0 \geq 0$ on $\Sigma (0)$. This means that $w^0 \equiv 0$, that is, $u(x) = u(-x_1 ,x' )$ for all $x= (x_1 ,x') \in \mathbb{R}^N $.   
Repeating again in every cartesian coordinate system centered at the origin, $u$ is radially symmetric and radially decreasing w.r.t. $0$. 
$\hfill \Box $
\section*{6. Examples of non-radial solutions in the case $q> q_S$}
\setcounter{section}{6}
\setcounter{equation}{0}
 
In this section we provide examples of non-radial solutions of problem {\bf (P)$_0$} when $N\geq 4$ and $q> q_S$. The first result concerns the Lane--Emden equation $-\Delta u = u^q $ and it has been obtained in  \cite{Danceretal}, Theorem 1.1.
\\[0.1cm]
{\bf Lemma 6.1. } 
{\sl Problem {\bf (P)}$_{\mathbf 0 }$ with $f(|x|,u)= u^q $ has infinitely many nonradial  solutions if $N\geq 4$ and 
\begin{equation}
\label{nonradial1}
q_S <q< \left\{ 
\begin{array}{ll}
+\infty & \mbox{ if $\ 4 \leq N\leq 11$}
\\
\frac{ (N-3)^2 -4N +4 +8\sqrt{N-2}}{(N-3)(N-11)} & \mbox{ if $\ N\geq 12 $}
\end{array}
\right. 
.
\end{equation}
These solutions take the form
\begin{equation}
\label{unonrad}
u(x) = |x|^{-\gamma} V(x|x|^{-1} ) ,
\end{equation}
where $\gamma = 2/(q-1)$ and $V$ is a non-constant solution of (\ref{sphericallaplace})
with $a= \gamma (N-2-\gamma )$. 
}
\\[0.1cm]
 Next,  we construct non-radial solutions for the H\'enon equation $-\Delta u= |x|^{-\ell} u^q $. We will use the following
\\[0.1cm]
{\bf Lemma 6.2. } {\sl (see \cite{BrezisLi}, Remark 2) 
\\
Assume that $N\geq 4$ and 
$q>q_S$. Then there is a number $\varepsilon_0 >0$, such that  problem  (\ref{sphericallaplace}) has a non-constant solution whenever}
\begin{equation}
\label{nonradial2}
0< \frac{N-1}{q-1} -\varepsilon_0 < a< \frac{N-1}{q-1} . 
\end{equation} 
\\[0.1cm]
  Using Lemma 6.2 we obtain  
\\[0.1cm]
{\bf Lemma 6.3. } {\sl Assume that $N\geq 4$ and $q>q_S $. Then, with the number $\varepsilon _0 $ of Lemma 6.2,  problem {\bf (P)}$_{\mathbf 0}$ with 
$$
f(|x|,u) = |x|^{-\ell} u^q 
$$
has a nonradial solution if 
\begin{eqnarray}
\label{nonradial3}
 & & 2-
\frac{(N-2)(q-1)}{2} +(q-1) 
\sqrt{ 
\left( \frac{N-2}{2}\right) ^2 - 
\frac{N-1}{q-1} 
}  
\\
\nonumber 
< \ell & < &      
2 -\frac{(N-2)(q-1)}{2} +(q-1) 
\sqrt{ 
\left( \frac{N-2}{2}\right) ^2 - 
\frac{N-1}{q-1}  + \varepsilon _0
}  
 .
\end{eqnarray} 
These solutions take the form (\ref{unonrad}) where
\begin{equation}
\label{defgamma}
\gamma  = \frac{2-\ell}{q-1}
\end{equation}
and $V$ is a non-constant solution of (\ref{sphericallaplace}) with 
\begin{equation}
\label{defa}  
a= \gamma (N-2 -\gamma ).
\end{equation}
}
{\bf Remark 6.4. } Due to our assumptions on $q$ and $\varepsilon $, (\ref{nonradial3}) implies that $\ell \in (0,2)$.
\\[0.1cm] 
{\sl Proof of Lemma 6.3 :} Let $\ell$ satisfy (\ref{nonradial3}), and define 
$
\gamma $  by (\ref{defgamma}) and $a$ by (\ref{defa})
Then, $a$ satisfies (\ref{nonradial2}) by our assumptions. Hence, 
Lemma 6.2 tells us that problem (\ref{sphericallaplace}) has a nonconstant solution for these values of $q$ and $a $. Furthermore, defining $u$ by (\ref{unonrad}), a short computation shows that $u$ is a solution of problem 
{\bf (P)}$_{\mathbf 0 }$ with $f(|x|, u) = |x|^{-\ell} u^q $.
$\hfill \Box $  
\section*{7. Appendix}
\setcounter{section}{7}
\setcounter{equation}{0}
 Assume that $\Omega \subset \mathbb{R} ^N $ is a bounded domain and $a_{ij}$, $b_i $ and $c$ are bounded 
measurable functions on $\Omega $ satisfying (\ref{ellcoeff}).
The following results are well known:
\\[0.1cm]
{\bf Theorem 7.1. }  {\sl (see \cite{BeNi}, p.4) $ \ $ Let $w\in W^{1,2} (\Omega )\cap 
W_{\mbox{loc}}^{2,N} (\Omega )$, and
\begin{equation}
\label{501}
-\sum_{i,j=1} ^N a_{ij} w_{x_i x_j } +\sum_{i=1}^N b_i w_{x_i } \leq cw \qquad \mbox{ in } \quad 
\Omega .
\end{equation}
{\bf (I)} Strong Maximum Principle :
If $w\leq 0$ in $\Omega $, $ x^0 \in \Omega $ and $w(x^0 )=0$ then $w\equiv 0$ in 
$\Omega $.
\\[0.1cm]
{\bf (II)}  Hopf's Boundary Point Lemma : 
If $w\leq 0$ in $\Omega $, $\partial \Omega$ is smooth in a neighborhood of $x^0 \in 
\partial \Omega ,\ w\in C^1 (\Omega \cup \{ x^0 \} )$ and $w(x^0 )=0$,
then
$$
\frac{\partial w}{\partial \nu } (x^0 ) \geq 0,\qquad \nu :\ \mbox{  
exterior normal },
$$
where equality is valid only if $w\equiv 0$ in $\Omega$.}
\\[0.2cm]
{\bf Acknowledgement:} This project was partially supported by South China University of Technology, Guangzhou (SCUT), by the Leverhulme Trust (UK), and project DAAD-ProMath 57425433. The authors want to thank SCUT, the Universidad de La Frontera (Temuco), and the Universit$\ddot{\rm a}$t Kassel for visiting appointments and for their kind hospitality.  


\small
\begin{thebibliography}{99}
\bibitem{AvilaBrock}
{\sc A. Avila, F. Brock}, Asymptotics at infinity of solutions for p-Laplace equations in
exterior domains. {\sl Nonlinear Analysis} {\bf 69} (2008), 1615--1628.
\bibitem{BeNi} 
{\sc H. Berestycki, L. Nirenberg}, 
On the method of moving planes and the sliding method.  
{\sl Bol.  Soc.  Bras.  Mat.} {\bf 22} (1991), 1--37.
\bibitem{BVV}
{\sc M.F. Bidaut-V\'eron, L. Veron}, 
Nonlinear elliptic equations on compact Riemannian manifolds and asymptotics of Emden equations. 
{\sl Invent. Math.} {\bf 106} (1991), 489--539. 
\bibitem{BidautPoho}
{\sc M.F. Bidaut-V\'eron, S. Pohozaev}, 
Nonexistence results and estimates for some nonlinear elliptic problems. 
{\sl J. Anal. Math. } {84} (2001), 1--49. 
\bibitem{BrezisLi}
{\sc H. Brezis, Yanyan Li}, 
Some nonlinear elliptic equations have only constant solutions. {\sl J. Partial Diff. Eqs.} {\bf 19} (2006), 208--217. 
\bibitem{CafGS}
{\sc L.A. Caffarelli, B. Gidas, J. Spruck}, Asymptotic symmetry and local behavior of semilinear elliptic equations with critical Sobolev growth. {\sl Comm. Pure Appl. Math. } {\bf 42} (1989), no. 3, 271--297.
\bibitem{cafflinirenberg}
{\sc L. Caffarelli, Y.Y. Li, L. Nirenberg}, Some remarks on singular solutions of nonlinear elliptic equations.
II: symmetry and monotonicity via moving planes, in: {\sl Advances in Geometric Analysis}, in: {\sl Adv. Lect.
Math. (ALM)}, vol. {\bf 21}, Int. Press, Somerville, MA, (2012), 97--105.
\bibitem{ChenLi}
{\sc Wenxiong Chen, Congming Li},
Methods on nonlinear elliptic equations. 
{\sl AIMS Series on Differential Equations \& Dynamical Systems, 4}. AIMS, Springfield, MO, (2010). 299 pp.
\bibitem{DamPaRa}
{\sc L. Damascelli, F. Pacella, M. Ramaswamy}, 
Symmetry of ground states of 
$p$-Laplace equations via the Moving Plane Method. 
{\sl Arch. Rat. Mech. Anal.} {\bf 148} (1999),  291--308.
\bibitem{DamRa} 
{\sc L. Damascelli, M. Ramaswamy}, 
Symmetry of $C^1$-solutions of $p$-Laplace equations in $\mathbb{R}^N$. 
{\sl Advances in Nonlinear Studies}, vol. 1, no. 1 (2001), 40--64.
\bibitem{DamSci}
{\sc L. Damascelli, B. Sciunzi}, 
Regularity, monotonicity and symmetry of positive solutions of $m$-Laplace equations.
{\sl J. Differential Equations} {\bf 206} (2004), 483--515.
\bibitem{DaDuGuo}
{\sc E.N. Dancer, Yihong Du, Zongming Guo}, 
Finite Morse index solutions of an elliptic equation with supercritical exponent. 
{\sl Journal Differential Equations} {\bf 250} (2011), 3281--3310. 
\bibitem{Danceretal}
{\sc E.N. Dancer, Zongming Guo, Juncheng Wei}, 
Non-radial singular solutions of Lane-Emden equations in $\mathbb{R}^N $. 
{\sl Indiana Univ. Math. J.} {\bf 61} (2012), no.3, 1971--1996.
\bibitem{Dengetal}
{\sc Yinbin Deng, Yi Li, Fen Yang},
On the positive radial solutions of a class of singular 
semilinear elliptic equations. 
{\sl Journal Differential Equations} 
{\bf 253} (2012), 481--501. 
\bibitem{Fra}
{\sc L.E. Fraenkel}, 
An introduction to maximum principles
and symmetry in elliptic problems.
{\sl Cambridge University Press } (2000), 1--351.
\bibitem{GNN}
{\sc B. Gidas, Wei-Ming Ni, L. Nirenberg},
Symmetry and related properties via the maximum principle.
{\sl Comm. Math. Phys.} {\bf 68} (1979), 209--243.
\bibitem{GNNrn}
{\sc B. Gidas, Wei-Ming Ni, L. Nirenberg},
Symmetry of positive solutions of nonlinear elliptic equations in $\mathbb{R}^N $,
in:
{\sl Mathematical Analysis and Applications, Part A},
in:
Adv. in Math., Suppl. Stud., vol 7a, Academic Press, New York, London, 1981, 369--402.
\bibitem{GidSpruck}
{\sc B. Gidas, J. Spruck},
Global and local behavior of positive solutions
of nonlinear elliptic equations.
{\sl Comm. Pure Appl. Math.} {\bf 34} (1981), 525--598.
\bibitem{Guo}
{\sc Zongming Guo}, On the symmetry of positive solutions of the Lane--Emden equation with supercritical exponent. {\sl Adv. Differential Equations} {\bf 7}, no.6 (2002), 641--666.
\bibitem{Hsia}
{\sc C.H. Hsia, C.S. Lin, Z.Q. Wang},
Asymptotic symmetry and local behaviors of solutions to a class of anisotropic elliptic equations. 
{\sl Indiana Univ. Math. J.} {\bf 60} (2011), no. 5, 1623--1654. 
\bibitem{CongmingLi}
{\sc Congming Li}, Monotonicity and symmetry of solutions of fully nonlinear elliptic equations on unbounded domains. {\sl Comm. Partial Differential Equations} {\bf 16} (1991), no. 4-5, 585--615.
\bibitem{YiLi}
{\sc Yi Li}, On the positive solutions of the Matukuma equation. {\sl Duke Math. J.} {\bf 70} (1993), no. 3, 575--589.
\bibitem{LiNi0}
{\sc Yi Li, Wei-Ming Ni}, On the existence and symmetry properties
of finite total mass solutions of the
Matukuma equation, the Eddington equation
and their generalizations.
{\sl Arch. Rational Mech. Anal.} {\bf 108} (1989), no. 2, 175--194. 
\bibitem{LiNi1}
{\sc Yi Li, Wei-Ming Ni}, On the asymptotic behavior and radial symmetry of positive solutions of semilinear elliptic equations in
$\mathbb{R}^N $. I. Asymptotic behavior,
{\sl Arch. Rat. Mech. Anal.} {\bf 118} (3) (1992), 195--222.
II. Radial symmetry.
{\sl Arch. Rat. Mech. Anal.} {\bf 118} (3) (1992), 223--243.
\bibitem{LiNi2}
{\sc Yi Li, Wei-Ming Ni}, Radial symmetry of positive solutions of nonlinear elliptic equations in $\mathbb{R} ^N $. {\sl Commun. Partial Differential Equations} {\bf 18} (1993), No.5--6, 1043--1054. 
\bibitem{MitPoh}
{\sc E. Mitidieri, S.J. Pohozaev},
A priori estimates and the absence of solutions of nonlinear partial differential equations and inequalities. (Russian) {\sl Tr. Mat. Inst. Steklova } {\bf 234} (2001), 1--384; translation in: {\sl Proc. Steklov Inst. Math. } {\bf 234} (2001), no. 3, 1--362.
\bibitem{Naito} 
{\sc Y. Naito}, 
Radial symmetry of positive solutions for semilinear elliptic equations in $\mathbb{R}^N$.
{\sl J. Korean Math. Soc.}  {\bf 37} (2000), no.5, 751--761.
\bibitem{weimingni}
{\sc Wei-Ming Ni}, 
Qualitative properties of solutions to elliptic problems. 
{\sl Handbook of Differential Equations;
Stationary Partial Differential Equations, Vol.1 } (2004),
157--233.
\bibitem{PhanSouplet}
{\sc Q. Phan, Ph. Souplet},
Liouville-type theorems and bounds of solutions of Hardy-H\'{e}non equations. 
{\sl J. Differential Equations} {\bf 252} (2012), 2544--2562.
\bibitem{PoQuSo}
{\sc P. Polacik, P. Quittner, Ph. Souplet}, Singularity and decay estimates in superlinear problems via Liouville-type theorems.
{\sl Duke Math. J.} {\bf 139}, no.3 (2007), 555--579.  
\bibitem{ProtWein}
{\sc M.H. Protter, H.F. Weinberger}, 
Maximum principles in differential equations. 
{\sl Springer-Verlag}, N.Y. (1984).
\bibitem{sciunzi}
{\sc B. Sciunzi}, 
On the moving plane method for singular solutions to semilinear elliptic equations. {\sl J. Math.
Pures Appl.} (9) {\bf 108} (2017), 111--123.
\bibitem{Se}
{\sc J. Serrin},
A symmetry problem in potential theory.
{\sl Arch. Rat. Mech. Anal.} {\bf 43} (1971), 304--318.
\bibitem{SeZ} 
{\sc J. Serrin, H. Zou}, 
Symmetry of ground states of quasilinear elliptic equations. 
{\sl Arch. Rat. Mech. Anal.} {\bf 148} (1999), 265--290.
\bibitem{terracini}
{\sc S. Terracini}, 
On positive entire solutions to a class of equations with a singular coefficient and critical
exponent, {\sl Adv. Differ. Equ.} {\bf 1} (1996), 241--264.
\bibitem{Z1}
{\sc H. Zou},
Symmetry of positive solutions of $\Delta u + u^p =0 $ in $\mathbb{R}^N $,
{\sl J. Differential Equations} {\bf 120} (1) (1995), 46--88.
\bibitem{Z2}
{\sc H. Zou},
Slow decay and the Harnack inequality for positive solutions of $\Delta u + u^p =0 $ in $\mathbb{R}^N $, 
{\sl Differential Integral Equations} {\bf 8} (6) (1995), 1355--1368.
\bibitem{Z3}
{\sc H. Zou},
Symmetry of ground states of semilinear elliptic equations with mixed Sobolev growth.
{\sl Indiana Univ. Math. J.} {\bf 45} (1) (1996), 221--240.
\end{thebibliography}
\end{document}